%% file: paper.tex
\documentclass[review,onefignum,onetabnum,onealgnum]{siamart190516}
\usepackage{amsfonts}
\usepackage{graphicx} 
\usepackage{epstopdf}
\usepackage[caption=false]{subfig}

\usepackage{algorithm}
\usepackage{algorithmicx}
\usepackage{algpseudocode}

\ifpdf
\hypersetup{
  colorlinks=true,
  linkcolor=blue,
  citecolor=blue,
  urlcolor=blue,
  pdftitle={The Polylogarithm Function in Julia},
  pdfauthor={M Roughan}
}
\fi

\def\equationautorefname~#1\null{%
  (#1)\null
}

\newlength{\figwidth}
\input{lecture_commands}
\newcommand{\Li}{\mbox{Li}}

\title{The polylogarithm function in Julia} % \\ in IEEE-754 64-bit floating point.
\author{Matthew Roughan\thanks{ARC Centre of Excellence for
    Mathematical \& Statistical Frontiers in the School of
    Mathematical Sciences at the University of Adelaide,
    Australia. \email{matthew.roughan@adelaide.edu.au}} }

\begin{document}

\maketitle
 
\input{abstract}


% https://mathscinet.ams.org/mathscinet/msc/pdfs/classifications2010.pdf

\input{introduction.tex}

\input{background.tex}

\input{related.tex}

\input{computation.tex}

\input{algorithm.tex}

\input{tests.tex}

\input{conclusion.tex}

\section*{Acknowledgements}

We would like to thank the Australian Research Council for funding
through the Centre of Excellence for Mathematical \& Statistical
Frontiers (ACEMS), and grant DP110103505. I would also like to thank
Andrew Feutrill for help in revising the paper.

% \ifCLASSOPTIONcaptionsoff
%   \newpage 
% \fi 
  
% \clearpage
{\small
\bibliographystyle{siamplain}
\bibliography{polylogarithm,books,computation}
}\par\leavevmode

\end{document}

%% file: lecture_commands.tex
\newcommand{\nc}{\newcommand}
\nc{\eg}{{\it e.g., }} 
\newcommand{\ie}{{\em i.e., }}

\nc\E{\ensuremath{\mathbb{E}}}
\nc\V{\ensuremath{\mathbb{V}}}
\nc{\bb}{{\text{\boldmath$b$}}}
\nc{\bs}{{\text{\boldmath$s$}}}
\nc{\eps}{\varepsilon}
\nc{\g}{{$G(N,L)$ }}
\nc{\pp}{($P$)}
\nc{\dd}{($D$)}
\nc{\pps}{($P$)$\,\,$}
\nc{\dds}{($D$)$\,\,$}
\nc{\ep}{{\varepsilon}}
\nc{\bu}{{\text{\boldmath$u$}}}
\nc{\bw}{{\text{\boldmath$w$}}}
\nc{\bq}{{\text{\boldmath$q$}}}
\nc{\be}{{\text{\boldmath$e$}}}
\nc{\bd}{{\text{\boldmath$d$}}}
\nc{\cd}{$c_{ij}^d$}
\nc{\bh}{{\text{\boldmath$h$}}}
\nc{\bv}{{\text{\boldmath$v$}}}
\nc{\bx}{{\text{\boldmath$x$}}}
\nc{\by}{{\text{\boldmath$y$}}}
\nc{\bc}{{\text{\boldmath$c$}}}
\nc{\bo}{{\text{\boldmath$o$}}}
\nc{\bz}{{\text{\boldmath$z$}}}
\nc{\bzero}{{\text{\boldmath$0$}}}
\nc{\R}{\mathbb R}
\nc{\N}{\mathbb N}
\nc{\Z}{\mathbb Z}
\nc{\Q}{\mathbb Q}
\nc{\C}{\mathbb C}
\nc{\mbP}{\mathbb P}
\nc{\I}{\mathbb I}
\nc{\disp}{\displaystyle}
\nc{\br}{{\text{\boldmath$r$}}}
\nc{\bpi}{{\text{\boldmath$\pi$}}}
\nc{\ba}{{\text{\boldmath$a$}}}
\nc{\bp}{{\text{\boldmath$p$}}}
\nc{\ff}{{\mathbf f}}
\nc{\var}{V\!ar}
\nc{\dsum}{\displaystyle \sum}
\nc{\dlim}{\displaystyle \lim}
\nc{\dmin}{\displaystyle \min}
\nc{\dsup}{\displaystyle \sup}
\nc{\dinf}{\displaystyle \inf}
\nc{\cf}{$C(\ff)$}
\nc{\hc}{$\nabla^2C(\ff)$}
\nc{\gc}{$\nabla C(\ff)$}
\nc{\ul}{\underline}
\nc{\zb}{$\bar{\bz}$}
\nc{\pc}{\frac{\partial C(\ff)}{\partial f_e}}
\nc{\xu}{ \bx_{\mu}}
\nc{\dx}{$\Delta \bx$}
\nc{\du}{$\Delta \bx_{\mu}$}
\nc{\p}[1]{\hbox{p} \! \left( #1 \right)}
\nc{\pr}[1]{\hbox{P} \! \left( #1 \right)}
\nc{\prc}[2]{\hbox{P} \! \left( \left. #1 \; \right| #2 \right)}
\nc{\ec}{\end{mathcal}}

\numberwithin{equation}{section}
\numberwithin{figure}{section}

\nc\urlb{\begingroup \urlstyle{tt}\Url}

\newcounter{myquote_counter}
\newsavebox{\myquotebox}
\newsavebox{\myquoterefbox}
\newlength{\quotewidth}
\newlength{\quoterefindent}
\newcommand*{\scaledquotefactor}{0.95}
{\renewcommand*{\scaledquotefactor}{1.0}%
   \begin{lrbox}{\myquoterefbox}\begin{minipage}[b]{0.8\quotewidth}%
       \em #1%
   \end{minipage}%
   \end{lrbox}%
   \begin{lrbox}{\myquotebox}\begin{minipage}[b]{1.0\quotewidth} 
}%
{\end{minipage}\end{lrbox}%
\begin{center}%
 \colorbox{quotebackcolor}{\begin{minipage}{\quotewidth}
                        \scalebox{\scaledquotefactor}{\usebox{\myquotebox}}\\[1pt]
          \hspace*{\quoterefindent}\scalebox{\scaledquotefactor}{\usebox{\myquoterefbox}}
       \end{minipage}
      }      
\end{center}%
\par}

%% file: abstract.tex
\begin{abstract}
The polylogarithm function is one of the constellation of important mathematical functions. It has a long history, and many connections to other special functions and series, and many applications, for instance in statistical physics. However, the practical aspects of its numerical evaluation have not received the type of comprehensive treatments lavished on its siblings. Only a handful of formal publications consider the evaluation of the function, and most focus on a specific domain and/or presume arbitrary precision arithmetic will be used. And very little of the literature contains any formal validation of numerical performance. In this paper we present an algorithm for calculating polylogarithms for both complex parameter and argument and evaluate it thoroughly in comparison to the arbitrary precision implementation in Mathematica. The implementation was created in a new scientific computing language Julia, which is ideal for the purpose, but also allows us to write the code in a simple, natural manner so as to make it easy to port the implementation to other such languages.
\end{abstract}

%% file: introduction.tex
\section{Introduction}

\begin{figure}[th]
  \centering
  \subfloat[Integer parameter $s$ and real argument $z$.]{\includegraphics[height=0.4\textwidth]{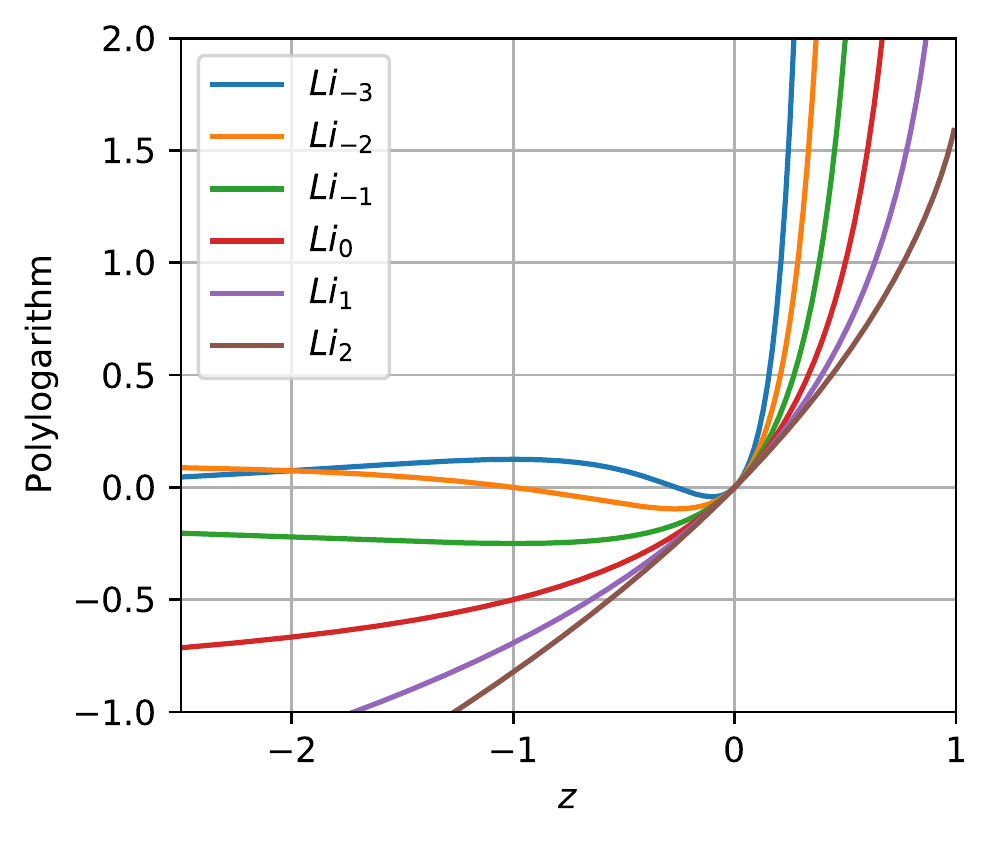}}
  \hspace*{3mm}
  \subfloat[Phase plot of $\Li_{-2}(z)$.]{\includegraphics[height=0.4\textwidth]{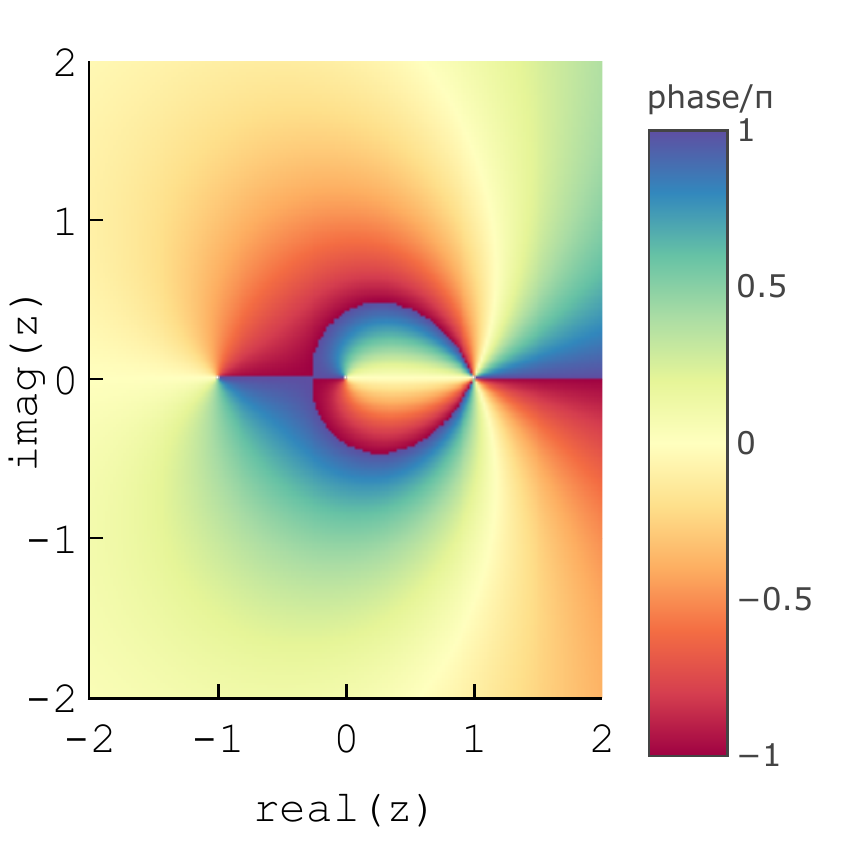}}
%  \subfloat[Phase of $\Li_{-2}(z)$.]{\includegraphics[height=0.4\textwidth]{Plots/polylog_ex_2.png}}
  % \begin{subfigure}[b]{width=0.49\textwidth}
  %   \centering
  %   \includegraphics[width=\textwidth]{Plots/polylog_ex_1.pdf}
  %   \caption{Integer parameter $s$ and real argument $z$.}
  %   \label{fig:polylog_ex_real}
  % \end{subfigure}
  \caption{Examples of the polylogarithm function.}
  \label{fig:polylog_ex}
\end{figure}

The polylogarithm function is defined by the sum
\begin{equation}
  \label{eq:def}
  \Li_s(z) =  \sum_{k=1}^{\infty} \frac{z^k}{k^s},
\end{equation}
for $|z| < 1$ (or for $|z| \leq 1$ when $\Re(s) \geq 2$), and by
analytic continuation to the entire complex plane. For instance
when $\Re(s)>0$ we can define it using the integral
\begin{equation}
  \label{eq:continuation}
  \Li_s(z) =  \frac{1}{z} \int_{0}^{\infty} \frac{t^{s-1}}{e^t/z - 1} dt,
\end{equation}
except for a pole at $z=1$ for $\Re(s) < 2$. There are many other 
representations for the function, but these suffice for our
understanding here. \autoref{fig:polylog_ex} shows some examples for
integer parameter $s$.

The polylogarithm is an interesting
function. Zagier~\cite{zagier07:_dia} describes just the dilogarithm
(the instance with $s=2$) as one of the simplest, and yet strangest of
non-elementary functions. He states
\begin{quote}
  ``almost all the formulas relating to it, have something of the
  fantastical in them, as if this function alone among all others
  possessed a sense of humor.''
\end{quote}
The polylogarithm has a long history; its early variants go back to
1696 with correspondance between Leibniz and the
Bernoullis\footnote{See Maximon~\cite{maximon03} for a brief
  history.}. And it is an important function. It has direct
relationships with the gamma function, Hurwitz and Riemann zeta
functions and many others. An entire book has been written on it
\cite{lewin81:_polyl}. It's relationships to Fermi-Dirac integrals
lead to physical applications, for instance, computations in
statistical mechanics \cite{lee97:_polyl_rieman}. It also has many
other applications including those in number theory and
geometry~\cite{maximon03}.

Note, however, that the polylogarithm function is unrelated to {\em
  polylogarithmic complexity} as discussed in complexity theory.
 
The function gains its name~\cite{zagier07:_dia} by comparison to the
Taylor series of the ordinary logarithm
\[ - \ln(1-z) = \sum_{n=1}^{\infty} \frac{z^n}{n}, 
\]
leading to some authors (\eg \cite{maximon03}) defining the
polylogarithm to refer to the integer cases $s=n$, and using the term
Jonqui\`{e}re's function or {\em fractional} polylogarithm for
non-integer cases. We use the term polylogarithm to include all values
of $s \in \C$.

There are many publications on the function, but only a handful on its
numerical
evaluation~\cite{jacobs72:_numer_calcul_polyl,wood92:_comput_poly,crandall06:_note,vepstas08:_hurwit,crandall12:_unified,BAILEY2015115,bailey15:_crandall}. Zagier~\cite{zagier07:_dia}
precisely describes the situation
\begin{quote}
  ``It occurs not quite often enough, and in not quite an important
  enough way, to be included in the Valhalla of the great
  transcendental functions---the gamma function, Bessel and
  Legendre-functions, hypergeometric series, or Riemann's zeta
  function. And yet it occurs too often, and in far too varied
  contexts, to be dismissed as a mere curiosity.''
\end{quote}
It occurs not quite often enough to have always been given the
consideration it deserves. In particular little work has been done on
its numerical evaluation, and all of the existing works are
incomplete, contain errors, or are specific to particular domains or
sets of parameters, \eg $s$ integer and real (see
\autoref{sec:related} for details).  Notably, almost none of these
present anything serious in the way of validation. And most existing
implementations are (i) not open, and/or (ii) use arbitrary-precision
arithmetic~\cite{polylog_python,polyl_wolfr_languag_system,vepstas19:_anant},
in order to produce high-precision results for instance for examining
zeros. Instead, we aim to provide a standard IEEE floating
point~\cite{ieee754} implementation of reasonable accuracy used for
instance for computing zeta-distribution moment-generating
functions. Obtaining reasonable accuracy over the whole parameter
space within the limits of 64 bit precision is still challenging. For
instance, some stopping rules that are relatively straight forward in
the deep tail of a series are less reliable early on in the sequence.

The implementation is provided in a relatively new programming
language -- Julia~\cite{bezanson17:_julia} -- designed specifically
for numerical computing. Julia is an ideal language for such tasks,
providing both an adaptable and dynamic high-level language, but also
very good computational
performance~\cite{karpinski14:_man_creat_one_progr_languag}.  Note,
however, that our code is open
source\footnote{Source available at \url{https://github.com/mroughan/Polylogarithms.jl}.} under
the MIT licence and the intent is that the code be straight-forward
enough to be adapted easily to any modern procedural language.

The algorithm is tuned to provide relative accuracy better than
$10^{-12}$ using standard IEEE double-precision floating point
calculations. In over 30,000 tests it fails to attain this accuracy in
only 88 cases, and the worst case accuracy is $1.1 \times
10^{-11}$. The typical accuracy is much better. The implementation's
speed is comparable or substantially better than the alternatives for
which we have data.

An additional goal of this work is to present a complete set of
information for implementation of this function in other languages and
settings, from the ground up, including many of the small technical
pieces that are missing from earlier expositions.

% The definition is sufficient for computation where it converges, but
% it does not for a large part of the complex plane, and its convergence
% leaves much to be desired in some domains. Apart from simple
% convergence, there are certain special points that require special
% care, in particular points $s$ near positive integers. At these
% points, simplistic summations may contain elements with poles, though
% they cancel leading to finite sums.

%% file: background.tex
\section{Notation, Conventions and Standard Results}
\label{sec:back}

The standard notation for the polylogarithm is $\Li_s(z)$, where we
refer to $s$ as the parameter, and $z$ as the input argument.  In the
case where $s$ is a real integer, it is often represented by $n$.

Note that the notation $\Li$ has also been used for the Eulerian
logarithmic integral, which is not under consideration here.

\subsection{Standard Functions}

There is a common set of useful functions -- Zagier's Valhalla --
defined in most computational packages. Most are well-known, \eg see
\cite{Abramowitz_and_Stegun}, but we define our notation here to be
precise.
\begin{eqnarray*}
  \Gamma(s) & = & \mbox{the gamma function}, \nonumber \\
  \psi(n)    & = & \mbox{the digamma function} = \mbox{the derivative of log-gamma},  \\
  \psi^{(m)}(n) & = & \mbox{the polygamma function} = \mbox{the $m$th derivative of log-gamma}, \\
  \zeta(s) & = & \mbox{the Riemann zeta function} =  \sum_{n=0}^{\infty} 1/n^s,
                  \mbox{ for }  \Re(s)>1,        \label{eq:zeta} \\
  \zeta(s,q) & = & \mbox{the Hurwitz zeta function} =  \sum_{n=0}^{\infty} 1/(n+q)^s,
                  \mbox{ for }  \Re(s)>1, \Re(q)>0,      \label{eq:hurwitz} \\
  \eta(s)  & = & \mbox{the Dirichlet eta function}  = \big( 1 - 2^{1-s} \big)\zeta(s),  \label{eq:eta} \\
  \beta(s) & = & \mbox{the Dirichlet beta function} =
                 4^{-s} \big( \zeta(s,1/4) - \zeta(s,3/4) \big). \label{eq:beta}
\end{eqnarray*}
These functions are provided by, or calculable directly from standard
numerical packages such as the {\tt
  SpecialFunctions}\footnote{\url{https://juliamath.github.io/SpecialFunctions.jl/}}
package in Julia.

% In many cases functions are extended outside the region of convergence
% by analytic continuation.

% Note the sign of $\eta(s)$ is different from Wood.

% http://mathworld.wolfram.com/PolygammaFunction.html
These have many useful relationships: \eg one that is
used here is 
\cite[6.4.2]{Abramowitz_and_Stegun}
\begin{equation}
  \label{eq:zeta_polygamma}
  \psi^{(n)}(1) = (-1)^{n+1} n! \,\zeta(n+1), \mbox{ for } n = 1,2,\ldots.
\end{equation}
% \begin{eqnarray}
%   \label{eq:psi}
%   \psi(n)    & = & \mbox{digamma function} \\
%   \psi^{(m)}(n)    & = & \mbox{$m$th derivative of digamma} \\
%                    & = & \mbox{the polygamma function} \\
% \end{eqnarray}

\subsection{Additional Functions and Sequences}

There are a number of sequences and functions used in computing
polylogarithms or as reference points for testing that are not as
commonly implemented in standard numerical packages (\eg these are not
provided by the {\tt SpecialFunctions} package in Julia).  In
particular:
\begin{itemize}
\item The {\em Stieltjes constants}~\cite{liang72:_stiel_const}
  $\gamma_N$, which we provide in our code from values from the
  OEIS~\cite{oeis}. Note that $\gamma_0 = \gamma$, which is also
  called the Euler-Mascheroni constant. The sequence is sometimes
  called the generalised Euler constants, but that this is confusing
  because there are other series called Euler numbers. Euler was
  prolific.
  
\item The {\em harmonic numbers}, ${\cal H}_n = \sum_{k=1}^{n} 1/k$ (where
  ${\cal H}_0=0$ by convention) and {\em generalized harmonic
    numbers}, 
  ${\cal H}_{n,r} = \sum_{k=1}^{n} 1/k^r$. We calculate values for
  smaller $n$ directly, while for larger we use the identity that
  \[ {\cal H}_n = \psi(n+1) + \gamma,
  \]
  and we use a similar identity for generalised harmonic numbers.

\item The {\em Bernoulli numbers} $B_n$ and polynomials $B_n(x)$. The
  first 35 numbers are provided as exact rationals\footnote{Julia
    provides a {\tt Rational} number type.}  based on
  \cite[A027642]{oeis}, and larger values are derived as real numbers
  from the polynomials using $B_n = B_n(0)$ because the integers in
  the rational representation come close to the bounds for 64 bit
  integers. Bernoulli polynomials are given exactly up to 5th order,
  otherwise calculated using the standard recursion, or for large
  $n$ using the following identity \cite{coffey09:_hurwit},
  $ B_n(x) = -n \zeta(1-n, x).  $

\end{itemize}
The advantage of expressing harmonic number and Bernoulli polynomials
in terms of other special functions is that these special functions
are given in many now standard numerical libraries, for instance, the
{\tt SpecialFunctions} module of Julia, but for small $n$ it can be
faster and more precise to calculate them directly.

\subsection{The Branch}
\label{sec:branch}

The polylogarithm function has a branch on the real axis for
$z \in [1,\infty)$. The conventional behaviour around the branch is
exemplified by particular values such as given in
\cite{crandall06:_note} for infinitesimally small $\epsilon > 0$
\begin{eqnarray}
  \label{eq:discontinuityL2}
  \Li_2(2)              & = & \frac{\pi^2}{4} - i \pi \ln 2, \\
  \Li_2(2 + i \epsilon) & = & \frac{\pi^2}{4} + i \pi \ln 2. 
\end{eqnarray}
In general the discontinuity should take the form
\begin{equation}
  \mbox{Disc} \; \Li_s(z) = 2 \pi i \frac{(\ln z)^{s-1}}{\Gamma(s)},
\end{equation}
the convention being \cite[(3.1)]{wood92:_comput_poly} that 
\begin{eqnarray}
  \label{eq:discontinuity}
   \Im\big( \Li_s(z) \big)              & = &  - \pi \frac{(\ln z)^{s-1}}{\Gamma(s)}, \\
   \Im\big( \Li_s(z + i \epsilon) \big) & = &  + \pi \frac{(\ln z)^{s-1}}{\Gamma(s)},
\end{eqnarray}
and these terms go to zero for $s = n \leq 0$.

% \subsection{Limits}

% Many limits of the polylogarithm function are known
%   \begin{eqnarray}
%     \label{eq:limit_z_0}
%     \lim_{z \rightarrow 0} \frac{\Li_s(z)}{z} & = & 1, \\
%    \label{eq:limit_}
%     \lim_{Re(\mu) \rightarrow \infty} \Li_s\left( e^{\mu} \right) & = & -\frac{\mu^s}{\Gamma(s+1)}
%     \mbox{ for } s \neq -1,-2, \ldots, \\
%    \label{eq:limit_}
%     \lim_{Re(\mu) \rightarrow \infty} \Li_{-n}\left( e^{\mu} \right) & = & -(-1)^n e^{-\mu} 
%     \mbox{ for } n = 1,2,3, \\
%     \label{eq:limit_s_infty}
%     \lim_{Re(s) \rightarrow \infty} \Li_s(z) & = & z. 
%   \end{eqnarray}
%   There are more known, but the central point is that the function
%   behaves relatively well in the limits of large $s$ and large $z$.

%% file: related.tex
\section{Related Work}
\label{sec:related}

It is surprising that there is relatively little written about
computing standard polylogarithms, much of it in informal literature,
rather than refereed publications. In rough chronological order:
\begin{itemize}

\item The pre-history of this work primarily concentrates only on di-
  and trilogarithms ($s=2$ and $3$), \eg see Jacobs and
  Lambert~\cite{jacobs72:_numer_calcul_polyl}.
  
\item Wood~\cite{wood92:_comput_poly} presents the first somewhat
  complete set of series for calculating the polylogarithm function in
  various domains. Apart from frequent typographical mistakes, there
  are several major difficulties with the manuscript.  Wood primarily
  presents series without any deep consideration of which should be
  applied.  Wood also focusses on $s$ real as do many other
  papers. The results often generalise, but special care was found to
  be needed for complex parameters.
  
  % The manuscript does not present an algorithm. This objection may
  % seem superficial, but there are subtle issues concerning which
  % series to use where, which many of the works on this calculation
  % ignore. Wood also focusses on $s$ real as do many other papers,
  % making glib statements that results extend to complex $s$ without
  % any verification.

\item Crandall~\cite{crandall06:_note} provides an actual algorithm,
  specifying choice of domain for each piece. However, in this paper
  Crandall only considers $s = n$, integer and real.  Crandall also
  provides only a few hints as to how the function is tested, but no
  numerical results. This is the first reference given in Python
  mpmath's implementation~\cite{polylog_python}, and is presumably the
  basis of the code there.
 
\item Vepstas~\cite{vepstas08:_hurwit} presents a new approach to
  solving problems of this type, but the algorithm uses
  arbitrary-precision arithmetic. We seek here to find a standard
  64-bit floating-point implementation. Vepstas provides more details
  of tests of the function than any of the other
  publications. However, even here, the results are often unclear: for
  instance, images without scales are published.
 
\item Crandall~\cite{crandall12:_unified} presents a unified version
  of his earlier approach for general $s$ (and other related
  functions), as well as the expression for the expansions near
  negative integer $s$, but in this work does not report a complete
  algorithm, and also suffers from some small mistakes. 

\item Bailey and Borwein~\cite{BAILEY2015115,bailey15:_crandall} use
  and refine Crandall's work for calculating polylogarithms and their
  derivatives. They discuss alternative approaches in different
  domains, but their main interest is in related functions and
  although they fill in some gaps of \cite{crandall12:_unified}, they
  do not present a complete algorithm either.
    
\end{itemize}

Existing software is also limited in this domain. Many implementations
are only valid for certain parameters (\eg integer and real $s$ or
even just a few values of $s$~\cite{voigt_polylogarith}).

Others
\cite{polylog_python,polyl_wolfr_languag_system,vepstas19:_anant} use
arbitrary-precision arithmetic, for instance in order to produce
high-precision results for computing zeros. Bailey and
Borwein~\cite{bailey08:_high_precis_comput_mathem_physic} argue
cogently for the need to have high-precision calculations in many
applications. However, a vast set of precedents show the value of
standard floating-point implementations (for instance see the C/C++
mathematics library) of useful functions. We here seek to add the
polylogarithm to the a list of commonly available functions. A simple
instance in which this would be useful is the calculation of the
moment generating and characteristic functions of the zeta
distribution.

% However, there does seem
% to be a need for simple, fast computations of typical floating point
% accuracy (since starting this project and making the initial code
% public, there have been several requests for it to be completed).

% We base our test comparisons on values derived from Mathematica's
% polylog function~\cite{polyl_wolfr_languag_system}.

% There are various generalisations and related topics that we do not
% consider here. 

%% file: computation.tex
\section{Components}

As with many special functions we use series representations to
calculate the polylogarithm. However, no single series converges over
the entire domain and there are some places where special numerical
care is needed even though the series technically converge. Thus, as
in previous works we will present several components. These are
largely consistent with those works, but some of the previous
works are incomplete or have mistakes, and hence we will list all of
the series used here in detail to eliminate any confusion.

We also use one identity to transform some parameter values into a
more amenable range.

In the following section we will precisely define how these are
combined into an algorithm. 

% The polylogarithm has many useful properties. We won't seek to list
% them all here, only those relevant for the computation and
% testing. For details and other properties see for instance
% \cite{lee97:_polyl_rieman,crandall06:_note,crandall12:_unified,}.

% Most importantly, it has a simple singularity at $z=1$ (and ), and a
% branch on the real axis for $z \in [1,\infty)$.

\subsection{Direct Series}

Our first port of call is the definition \autoref{eq:def}. This
definition can be used directly to calculate the polylogarithm for
$|z|<1$, and we shall refer to this approach to calculation as {\bf
  Series~1}. 

% \begin{figure}[ht]
%   \begin{center}
%     \includegraphics[width=\columnwidth]{Figs/domains.pdf}
%   \end{center}
%   \caption{A cartoon (the diagram is not to scale) of the domains in
%     which each series converges, and the partitioning stragies
%     suggested by Crandall~\cite{crandall06:_note} $a=1/2$, and Bailey
%     and Borwein~\cite{bailey15:_crandall} $a=1/4$. The lines show
%     regions of convergence; the shaded regions show which series is
%     used.}
%   \label{fig:domains}
% \end{figure}

Although it is convergent for $|z|<1$, the rate of convergence can be
quite slow close to the boundary.  Crandall~\cite{crandall06:_note}
suggests use of this series for $|z|<1/2$, which leads to an easy
calculation of 1 bit of precision per term (asymptotically) in the
series, based on the asymptotic dependence on the $|z|^k$ term.
Bailey and Borwein~\cite{bailey15:_crandall} suggest $|z|<1/4$ based
on experimental results, however we shall refine this the section
below. However, note that much of the commentary presumes very-high
precision will be required, and hence the sum must proceed deep into
the tail of the series. We are concerned here with finite precision
and when $\Re(s) < 0$ the early part of the sequence is strongly
impacted by the denominators of the sequence. If $z$ is small, we may
not need to proceed further into the tail, and so Crandall's comment
about bits per term may not apply.

% \subsection{}

% There is an expression for the polylogarithm in terms of the Hurwitz
% zeta function. Given Julia already implements this function it might
% be an obvious approach, but (i) this performed poorly in terms of
% accuracy, (ii) Julia's Hurwitz zeta function domain is limited; and
% (iii) the relationship fails near the points $s$ non-negative
% integers. Specific expansions are needed around these points, and the
% following alternative provides more immediate access to those
% expansions. Thus we will not use this approach.

\subsection{Alternative Series 2}
\label{sec:series2}

The main alternative power series is about $z=1$, where we
find~\cite[(32)]{crandall12:_unified} and
\cite[(9.3)]{wood92:_comput_poly} and \cite[(2)]{bailey15:_crandall}
\begin{equation}
  \label{eq:power_series2a}
  \Li_s(z) =  \Gamma(1-s) (-\ln z)^{s-1} 
              + \sum_{k=0}^{\infty} \frac{ \zeta(s - k) }{k!} (\ln z)^k,
\end{equation}
which converges for $|\ln z| < 2 \pi$ and $s$ not a positive integer. 
We call this {\bf Series 2}.

\begin{figure}[tb]
  \centering
  \includegraphics[width=0.49\columnwidth]{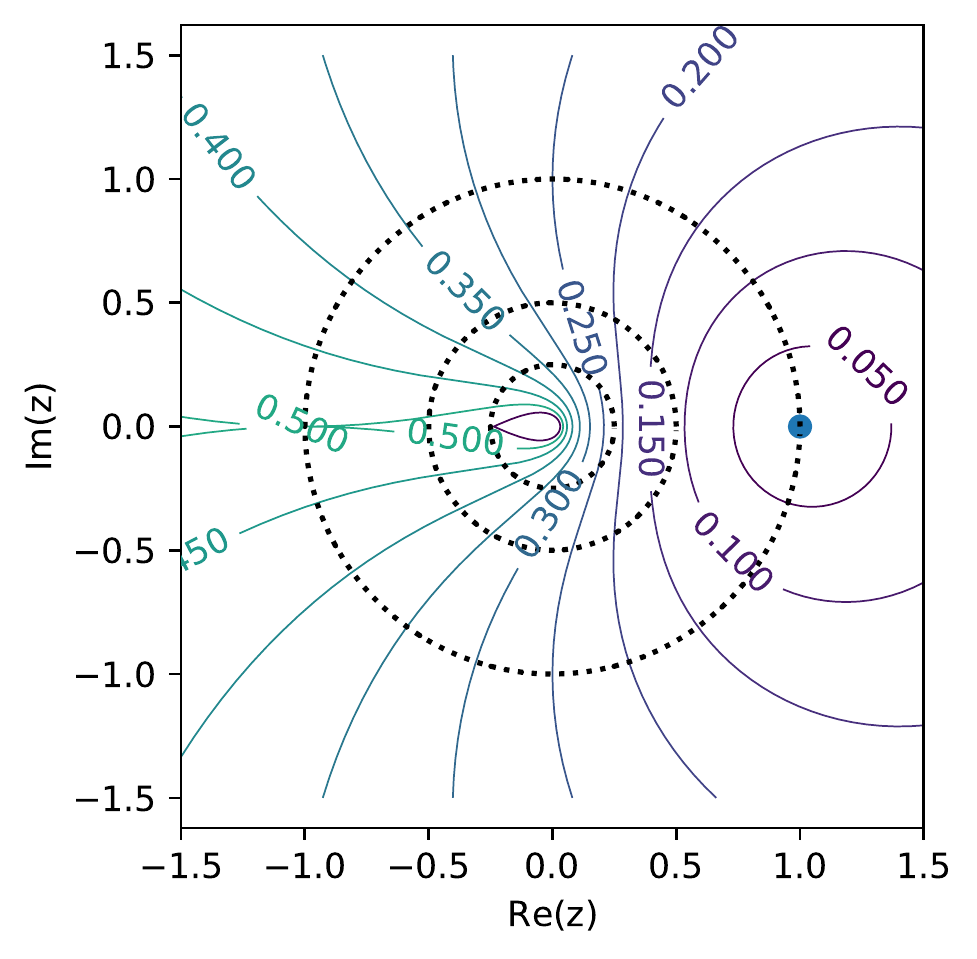}
  \includegraphics[width=0.49\columnwidth]{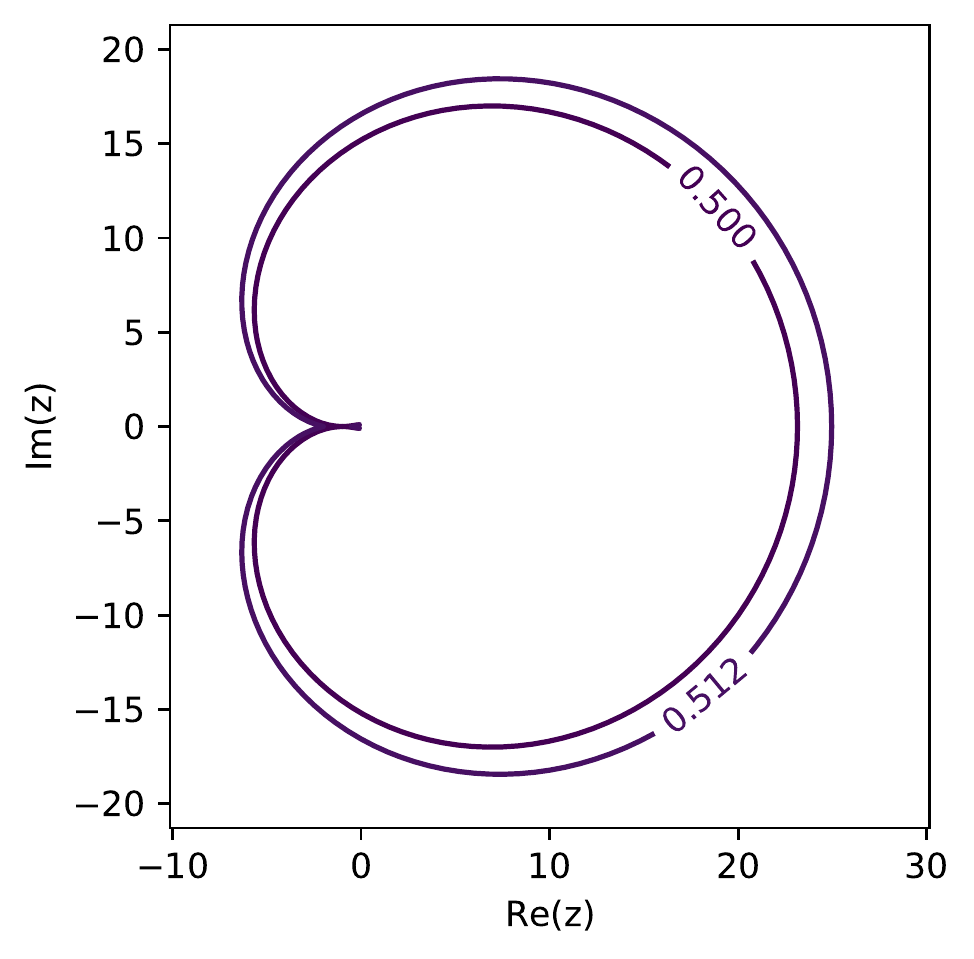}
   \caption{Contours of constant magnitude for $\ln z/2\pi$ (solid
    lines). The final tear-drop shaped contour near the origin
    corresponds to $\ln z/2\pi=0.55$. Dotted lines show
    $|z| = 0.25, 0.5$ and $1.0$ for comparison. The right figure
    illustrates the large regions of convergence outside the uni
    circle for given thresholds.} 
  \label{fig:logz_contours}
\end{figure}

Convergence in the region $|\ln z| < 2 \pi$ arises from Riemann's
functional equation \cite[23.2.6]{Abramowitz_and_Stegun}
\[ \zeta(s) = 2^s \pi^{s-1} \sin( \pi s/2) \Gamma(1-s) \zeta(1-s), \]
which, when substituted in the above series leads to
\begin{eqnarray*}
  \Li_s(z) & =  & \Gamma(1-s) (-\ln z)^{s-1}  \\
           & &   + 2^s \pi^{s-1} \sum_{k=0}^{\infty} \sin( \pi (s-k)/2) \zeta(1 - s + k)  \frac{\Gamma(1 - s + k) }{\Gamma(1+k)}
              \left(\frac{\ln z}{2 \pi} \right)^k,
\end{eqnarray*}
The terms $\zeta(1 - s + k)$ converges to 1 as $k \rightarrow \infty$,
and the sine is bounded, and hence the tail of the sequence has much
in common with that of Series~1, except that it is dominated by the
term $(\ln z/2\pi)^k$. \autoref{fig:logz_contours} shows contours of
$\ln z/2\pi$ both at small and large scales.
 
This series has a much larger domain of convergence than Series 1, but
there is a small region about the origin which is excluded and hence
we cannot discount Series~1.  Past approaches partition the space
using Series~1 inside some radius around the origin, and Series~2 in a
region around that with alternative radii being proposed by
Crandall~\cite{crandall06:_note} ($|z| = 1/2$) and Bailey and
Borwein~\cite{bailey15:_crandall} ($|z| = 1/4$).

However, it seems more natural (and is actually consistent with the
detailed results in \cite{bailey15:_crandall}) to make the dividing
line between Series~1 and Series~2 slightly more complex: we will
prefer Series~1 when $2 \pi |z| < | \ln z |$. This leads to a more
complex boundary, but results in a better tradeoff between the two
approaches in some regions, particularly on the negative real axis
where convergence of Series~2 is at its worst.

\begin{figure}[th]
  \centering
  \includegraphics[width=0.8\columnwidth]{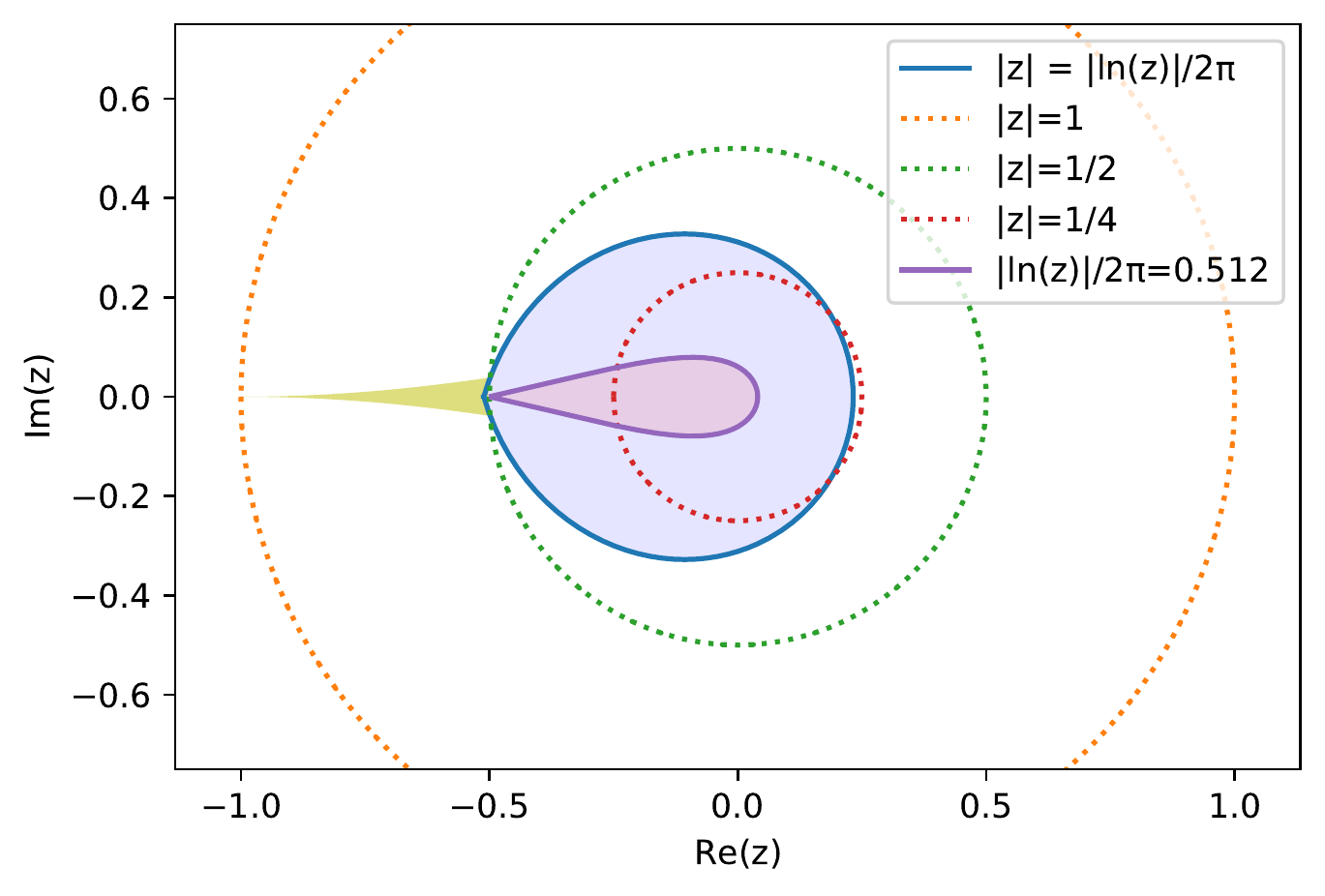}
  \caption{Boundaries in the choice of Series~1 and 2. Series~1 is
    used in the shaded blue region. The yellow region indicates the
    region that is missed if we bound use of Series~2 by
    $\ln z/2\pi \leq 0.5$. The shaded red tear-drop region indicates
    the region excluded if we use Series~2 when
    $\ln z/2\pi \leq 0.512$. As that this lies inside the region
    where Series~1 is used all of the unit disc is covered by
    either Series~1 or 2. }
  \label{fig:series_boundary}
\end{figure}

The regions are shown in \autoref{fig:series_boundary}. Note that on
the positive real axis we fall back to the advice of Bailey and
Borwein in that our transition occurs at $|z| \simeq 0.2323$ (near
1/4) and on the negative real axis we chose a threshold
$z \simeq 0.5113$ near that of Crandall (near 1/2). The shape of the
region proposed for the choice matches very closely to the empirical
results shown in Bailey and
Borwein~\cite[Figure~1]{bailey15:_crandall}.

Many existing works presume that use of Series~2 will be restricted to
the region $\ln z/2\pi \leq 0.5$ in order that, asymptotically, we
obtain 1 bit per term as for Series~1.
Figures~\ref{fig:logz_contours} and \ref{fig:series_boundary} together
illustrate a seemingly unstated fact that despite being a very large
region, the combined regions $|z| \leq 0.5$ and
$|\ln z|/2\pi \leq 0.5$ still leave a small gap in the unit disc near
the negative real axis (shaded yellow). This gap closes if we allow as
small a change as $|\ln z|/2\pi \leq 0.512$, and so we can either allow
a small increase in the region allowed for Series~2, or use an
alternative in this region. The duplication formula appears to work,
but we have not proved that it converges for every point in this
domain and so we prefer the former approach.

Detailed results of the tradeoffs are included in
\autoref{sec:test:trade12}.
  
% Given these thresholds we might expect performance to be worst along
% the negative real axis, but remember that in this domain signs will
% alternate resulting in some cancellation, so convergence may not
% suffer as much

\subsection{Alternative Series 3}

For $s=n>0$ (\ie a positive integer) we get the same summation,
however, the terms $\zeta(1)$ and $\Gamma(1-s)$ both have poles. The
poles cancel, but if we naively calculate these, we run into numerical
problems. Cancelling the two poles we get (as
in~\cite[(33)]{crandall12:_unified})
\begin{equation}
  \label{eq:power_series2x}
  \Li_{n}(z) =  \frac{(\ln z)^{n-1}}{(n-1)!} \big[ H_{n-1} - {\cal L}\big]
          + \sum_{k=0 \atop k\neq n-1}^{\infty} \zeta(n - k) \frac{(\ln z)^k}{k!},
\end{equation}
where ${\cal L} =\ln(-\ln z)$. This series could be used for cases
where $s=n$ but when $s$ is near an integer then there are numerical
difficulties in computing the difference of the two large terms. 
Hence for $s=n+\tau$ for small $\tau$, we perform an expansion around
$s=n$ to get
\begin{eqnarray}
  \label{eq:power_series2c}
  \Li_{n+\tau}(z) & = & \frac{(\ln z)^{n-1}}{(n-1)!} Q_{n-1}({\cal L},\tau) 
          + \sum_{k=0 \atop k\neq n-1}^{\infty}  \zeta(n+\tau - k) \frac{(\ln z)^k}{k!},
\end{eqnarray}
where the two problematic terms are grouped into $Q_{n}({\cal
  L},\tau)$, \ie 
\begin{eqnarray*}
  Q_{n}({\cal L},\tau) & = & \zeta(1) + (-1)^{n} n! \Gamma(-n-\tau) (\ln z)^\tau. 
\end{eqnarray*}
Crandall~\cite[(51)]{crandall12:_unified} expands $Q_{n}({\cal L},\tau)$
as a Taylor series about $\tau=0$ as 
\begin{eqnarray*}
  Q_{n}({\cal L},\tau)  & = & \sum_{j=0}^{\infty} c_{n,j}({\cal L}) \tau^j ,
\end{eqnarray*}
and gives the $c_{k,j}$ in \cite[pp.35-36]{crandall12:_unified}
recursively, but note that Crandall's manuscript has typographic
errors;  Bailey and Borwein~\cite{BAILEY2015115,bailey15:_crandall}
give the correct formula.
%
% by 
% \[ \label{eq:c_k}
%   c_{n,j}({\cal L}) = \frac{(-1)^j}{j!} \gamma_j - b_{k,j+1}({\cal L}) ,
% \]
% where
% \[
%   b_{k,j+1}({\cal L}) = \sum_{p+t+q=j \atop p,t,q \geq 0}
%                             \frac{{\cal L}^p}{p!}
%                             \frac{\Gamma^{(t)}(1)}{t!} (-1)^{t+q} f_{k,q}, 
% \]
% where $f_{k,q}$ is the calculated recursively using $f_{k,0}=1$ and 
% \[
%   f_{k,q} = \sum_{h=0}^{q} \frac{(-1)^h}{k^h} f_{k-1,q-h},
% \]
% and the derivatives of gamma functions at 1 are
% \begin{eqnarray*}
%   \Gamma^{(0)}(1) & = & 1, \\
%   \Gamma^{(1)}(1) & = & -\gamma, \\
%   \Gamma^{(2)}(1) & = & \gamma^2 - \gamma + \frac{\pi^2}{6}, \\
%   \Gamma^{(3)}(1) & = & -2 \gamma^3 + 9 \gamma^2 - (\pi^2 + 6) \gamma
%                         + \frac{3}{2} \pi^2 - 4 \zeta(3). 
% \end{eqnarray*}
However, for $|\tau| \ll 1$ we need only take a small number of terms,
the first three of which can be written explicitly as
\begin{eqnarray*}
  c_{n,0}({\cal L})  & = & H_n - {\cal L}, \\ 
  c_{n,1}({\cal L})  & = &
         -\gamma_1 - \frac{(\psi(n+1) -  {\cal L})^2}{2} - \left(\frac{\pi^2}{6} - \frac{\psi^{(1)}(n+1)}{2}\right), \\
  c_{n,2}({\cal L})  & = &
           \frac{\gamma_2}{2} + \frac{(\psi(n+1) -  {\cal L})^3}{6}\\
           &&     + (\psi(n+1) - {\cal L}) \left(\frac{\pi^2}{6} - \frac{\psi^{(1)}(n+1)}{2} \right)
                  + \frac{\psi^{(2)}(n+1)}{6}.
%       return stieltjes(2)   + d2^3/6 + d2*( π^2/6 - SpecialFunctions.polygamma(1,n+1)/2) + SpecialFunctions.polygamma(2,n+1)/6
\end{eqnarray*}
Noting that from \cite[6.3.2]{Abramowitz_and_Stegun} and
\autoref{eq:zeta_polygamma} we get relationships such as
\[ \psi^{(0)}(n) = \psi(n) = - \gamma + H_{n-1}
  \quad \mbox{ and } \quad 
    \psi^{(1)}(1) = \zeta(2) = \frac{\pi^2}{6},
\]
we see that $c_{n,0}({\cal L})$ and $c_{0,1}({\cal L})$ are consistent
with the formulas given in
\cite[pp.36]{crandall12:_unified}. Wood~\cite[(9.4)]{wood92:_comput_poly}
presents a similar expansion, but only the 1st term is the same. The
2nd term differs only in the signs of some terms and there are other
typographic errors in the work.
% These problems are avoided in the
% many implementations for $s=n$.

% suggesting a
% typo in those results\footnote{There is another typographic error in
%   the exponent of $(-w)$ in this equation as well, supporting that
%   conclusion.}.

% Q3 =  gen_euler(2) + (digamma(n)-lmu)^3/6 + (pi^2/6 + polygamma(1,n)/2)*(digamma(n)-lmu) + polygamma(2,n)/2

% , but for floating point calculations the
% first three suffice
% \begin{eqnarray}
%   c_{n,0}({\cal L}) & = & H_n - {\cal L}, \\
%   c_{0,1}({\cal L}) & = & -\gamma_1 - \frac{1}{2} \gamma^2 -
%                           \frac{\pi^2}{12} - \gamma {\cal L} -
%                           \frac{1}{2} {\cal L}^2,
% %   c_{n,2}({\cal L}) & = & 
% \end{eqnarray}
% The third term is not provided in closed form in
% \cite{crandall12:_unified}, but it is 

Taking $\tau=0$ we use just the first term leading to a result
consistent with \autoref{eq:power_series2x}. Taking a $j$th order approximatio we get
\[ 
    Q^{(j)}_n({\cal L},\tau)= c_{n,0}({\cal L}) + \tau c_{n,1}({\cal
      L}) + \cdots + \tau^j c_{n,j}({\cal L}),
\]
in 
\begin{equation}
  \label{eq:series3}
  \Li_{n+\tau}(z) \simeq \frac{(\ln z)^{n-1}}{(n-1)!} Q^{(j)}_{n-1}({\cal L},\tau)
          + \sum_{k=0 \atop k\neq n-1}^{\infty}  \zeta(n + \tau - k) \frac{(\ln z)^k}{k!}.
\end{equation}

\begin{figure}[th]
  \centering
  \includegraphics[width=0.9\columnwidth]{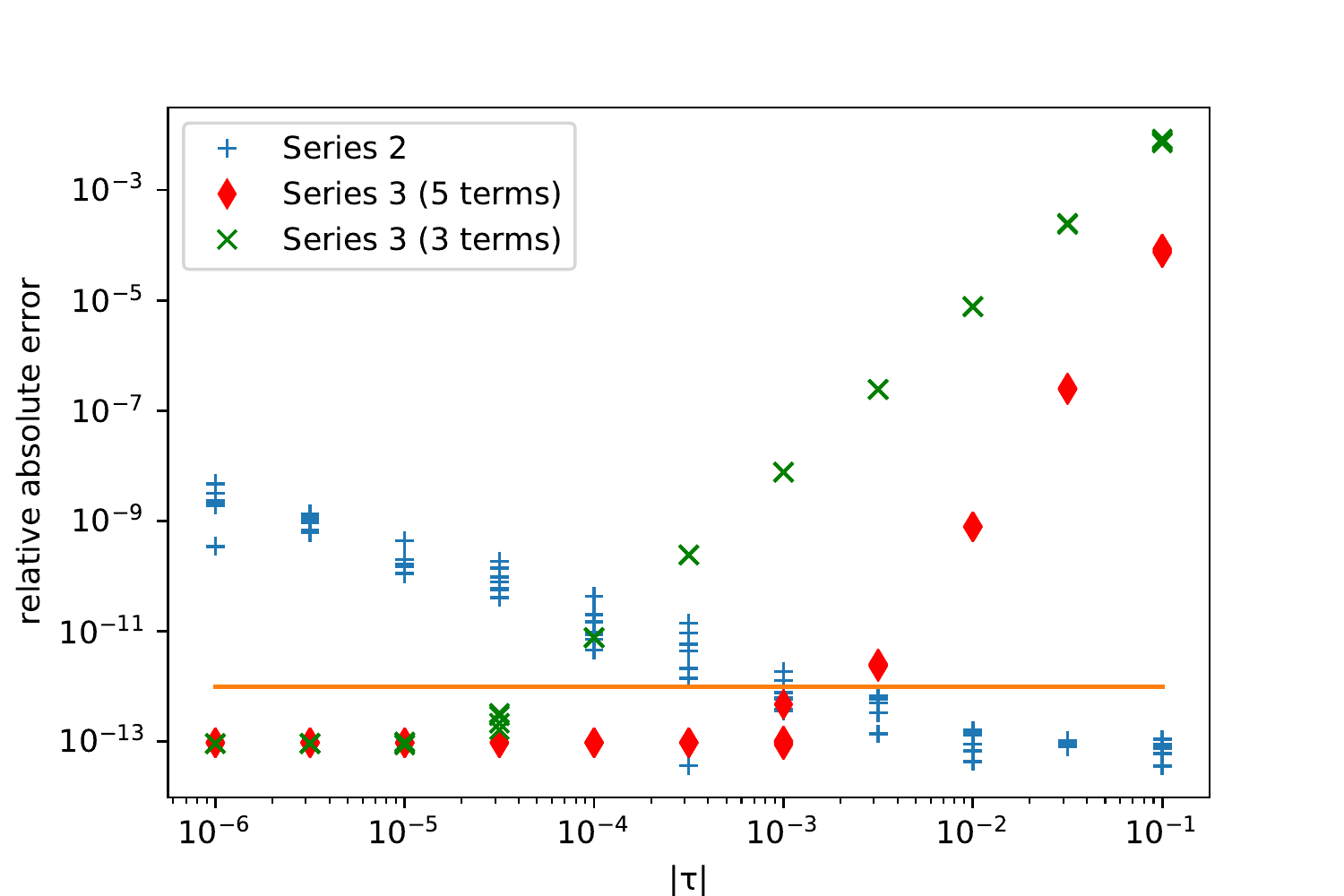}
  \caption{A comparison of Series~2 and 3 near positive integer
    parameters, \ie $s = n + \tau$, for $n>0$ (here $n=1$ and
    $z = -1/2$). The number of terms in Series~3 is with respect to
    the approximation of $Q_{n-1}({\cal L},\tau)$ used in
    \autoref{eq:power_series2c}. Note that the cross-over point
    between Series~2 and Series~3 achieves the desired precision for 5
    terms and a threshold around $10^{-3}$. The three term version of
    Series~3 is included because in that case the terms are known in
    closed form making calculation potentially faster with a loss of
    precision of about 1 order of magnitude. }
  \label{fig:convergence_b} 
\end{figure}

We refer to \autoref{eq:series3} as {\bf Series 3}. To determine how
many terms are needed in $Q$ and how small $|\tau|$ need be before we
swap to Series~3, we test the two series empirically near the point
$s=1$.  \autoref{fig:convergence_b} shows the errors as a function of
$|\tau|$.  We see quite similar results for several other values of
$s=n$ that were tested.  The figure indicates that around 5 terms are
needed before the cross-over point between the two series lies below
$10^{-12}$ (our target precision) and that the threshold should be
$\tau < 10^{-3}$.

\subsection{The Special Case At $z=1$}

There are a number of special cases where the polylogarithm resolves
to a more familiar function. The most obvious (the reason for which
the function is named) is that $\Li_1(z) = - \ln(1-z)$, but at the
moment we do not use special cases in the calculations so that we can
use these identies as test cases. However the following identity
\begin{equation}
  \label{eq:special_case_z=1}
  \Li_s(1) = \zeta(s), \; \mbox{ for } \Re(s) > 1, 
\end{equation}
is useful because at this allows Series~2 and Series~3 to be written
without considering the special case in detail. Note that when
$ \Re(s) \leq 1$ there is a pole at $z=1$, so the function returns
{\tt Inf}.
 
\subsection{Duplication Identity}

The polylogarithm satisfies the duplication identity~\cite[(e)]{lee97:_polyl_rieman}
\begin{equation}
  \label{eq:square}
  \Li_s(z) + \Li_s(-z) = 2^{1-s} \Li_s( z^2 ).
\end{equation}
Wood refers to this as the {\em square formula}
\cite[14.1]{wood92:_comput_poly}, and notes that in an alternative
\cite[15.1]{wood92:_comput_poly} form
\begin{equation}
  \label{eq:duplication}
  \Li_s(z) = 2^{s-1} \big(  \Li_s( \sqrt{z} )  + \Li_s( -\sqrt{z} ) \big],
\end{equation}
it could be use recursively to put $z$ into a domain where another
algorithm can be brought to bear. However Wood recommends against this
citing difficulties in implementing recursion in his contemporary Fortran
as well as the increase in the number of polylogarithms to be
calculated. Recursion is not difficult in many modern languages so
this can be a useful technique noting that the increase in the number
of polylogarithms is countered by faster convergence on the square
root.

% In particular if we consider Series~2 and 3, their convergence lies in
% the value of $a = \log|z|/ 2\pi$, and so although the applying
% \autoref{eq:square} once doubles the number of polylogarithms to be
% calculated, it reduces $a$ by a factor of two, and this is useful in
% the range, for instance $1/2 < \log|z|/ 2\pi 1$, where it will lead to
% much faster convergence, or in the range
% $1/4 < \log|z|/ 2\pi \leq 1/2$ where the rate of convergence is
% approximately doubled. This is quite helpful because the follow
% identity (which we resort to outside of these regions) is numerically
% finicky. 

\subsection{Reciprocal Identity}

The polylogarithm also satisfies another identity often called
Jonqui\`{e}re's Identity, which leads to a reciprocal relationship:
% that apply for all $s \in \C$
% \begin{eqnarray} 
% \lefteqn{\Li_s(z) + (-1)^s \Li_s(1/z)} \nonumber \\
% & = & \frac{(2 \pi i)^s}{\Gamma(s)} \zeta\left(1-s, \frac{1}{2} +  \frac{\ln(-z)}{2 \pi i}    \right) \nonumber  \\
% &   & - 2 \pi i  \Theta(z) \frac{(\ln z)^{s-1}}{\Gamma(s)},
%   \label{eq:reciprocal}
% \end{eqnarray}
\begin{eqnarray} 
\lefteqn{\Li_s(z) + (-1)^s \Li_s(1/z)} \nonumber \\
& = & \frac{(2 \pi i)^s}{\Gamma(s)} \zeta\left(1-s, \frac{1}{2} +  \frac{\ln(-z)}{2 \pi i}    \right) 
 - 2 \pi i  \Theta(z) \frac{(\ln z)^{s-1}}{\Gamma(s)},
  \label{eq:reciprocal}
\end{eqnarray}
where  $\Theta(z)$ is Crandall's {\em domain-dependent step function}
\begin{equation}
  \Theta(z) = \left\{
    \begin{array}{ll}
      1 & \mbox{ if } \Im(z) < 0 \mbox{ or } z \in [1,\infty), \\ 
      0 & \mbox{ otherwise}.
    \end{array}
    \right.
\end{equation}
The function $\Theta(z)$ is intended to provide the conventional
behaviour on and around the branch.

Wood~\cite{wood92:_comput_poly} and others suggest use of this for
large $|z|$. We tested the reciprocal computation but we found this
approach to be numerically unstable for $s$ with a large, negative
imaginary component\footnote{Wood and others primarily consider the
  calculation of the polylogarithm with real arguments.}  and thus
prefer not to use it here. This choice also reduces the dependence of
our code on the Hurwitz zeta function, which is perhaps the least
commonly implemented function (of those used here) in standard
scientific packages.

Note that if $s$ takes negative integer values we hit poles of
the Gamma functions. In these cases we use the relation
\cite[(10.3)]{wood92:_comput_poly}
\begin{equation} 
\Li_{-n}(z) + (-1)^n \Li_{-n}(1/z) = 0,
  \label{eq:reciprocal_neg_int}
\end{equation}
but we shall use this relationship primarily for validation in what
follows.

% In these cases we use the
% integer form given in \cite[1.3]{crandall06:_note}
% \begin{eqnarray} 
% \lefteqn{\Li_n(z) + (-1)^n \Li_n(1/z)} \nonumber \\
% & = & - \frac{(2 \pi i)^n}{n!} B_n\left(\frac{\ln z}{2 \pi i}    \right) 
%  - 2 \pi i  \Theta(z) \frac{(\ln z)^{n-1}}{(n-1)!},
%   \label{eq:reciprocal_int}
% \end{eqnarray}

% The boundary used to determine when this identity will
% be applied is determined by the optimal tradeoff between using
% Series~2 or Series~1 after applying the reciprocal (ignoring the cost
% of calculating the Hurwitz zeta) and thus is the reciprocal of the
% boundary computed above.

% When $s = n$ (a real integer) we can exploit
% \[ \zeta(-n,x) = -B_{n+1}(x)/(n+1), \] 
% for $Re(x)>0$ to obtain \mbox{\cite[(35)]{crandall12:_unified}}
% \begin{eqnarray} 
% \lefteqn{\Li_n(z) + (-1)^n \Li_n(1/z)} \nonumber \\
% & = & - \frac{(2 \pi i)^n}{n!} B_n\left( \frac{\ln z}{2 \pi i} \right)
%     -2  \pi i  \Theta(z) \frac{(\ln z)^{n-1}}{(n-1)!},
%   \label{eq:reciprocal_int}
% \end{eqnarray}
% for $z \not\in (0,1]$.

% As before we don't use exact integer relations such as this in
% computation, but will use it for verification.

%% file: algorithm.tex
\section{Algorithm}

\subsection{Dependencies}

The calculations require a number of other mathematical constants,
sequences and functions as described in \autoref{sec:back}. Many of
these (the gamma, digamma, polygamma, and Riemann, Hurwitz and
Dirichlet zeta functions) are provided by Julia's {\tt
  SpecialFunctions} package (v0.10.3) designed around the {\tt
  OpenSpecFun}\footnote{{\tt OpenSpecFun} uses AMOS and Faddeeva to
  provide Bessel, Airy and error functions.} and {\tt
  OpenLibm}\footnote{{\tt OpenLibm} uses the standard C libm, which
  includes, for instance the gamma function suite.} libraries.  The
implementation of the polygamma is restricted to integer orders $m$
but that is all we require. The Euler-Mascheroni constant $\gamma$ and
$\pi$ are provided as part of {\tt Base.MathConstants} as {\tt Irrational}
type constants, which is a type that supports both floating point and
functional (arbitrary precision) definition.

% Julia's comments says that it obtains the Hurwitz zeta function from
% the polygamma using the relationship \autoref{eq:zeta_polygamma},
% which the documentation notes that this is valid for $\Re z >
% 0$. However, the code itself uses direct evaluation of the defining
% power series over some domains and instead calculates the polygamma
% using this relationship. 

  % \[ \zeta(s, z) = \sum_{k=0}^\infty ((k+z)^2)^{-s/2} \]
  % where terms $k+z=0$ are excluded
  
Other required constants, sequences and functions are defined as part
of this package using standard values and algorithms as indicated in
\autoref{sec:back}. The implementations of these components are
comparatively straight forward.

% including
% Bernoulli sequences and polynomials, harmonic numbers, Stieltjes
% numbers and the Dirichlet beta function

\subsection{Domain breakup}

The algorithm selects various approaches in different domains.
The completed algorithm is given via the pseudo-code described in
\autoref{alg:main}, which describes the breakup of the input domain in
detail. The main breaks occur to use Series~3 near positive integer
values of $s$ and to separate the domains where Series~2 and Series~3
have the best convergence. 

There are special cases of the function for certain values (for
instance $\Li_1(z) = -\ln(1-z)$) but we do not use these at present
because these identities are useful in testing series convergence.
These might be used more in future versions.

% The main distinctions are
% \begin{itemize}
  
% \item Certain series converge faster within particular domains of $z$;
%   \autoref{fig:series_boundary} shows the breakdown.

% \item Series~2 behaves poorly near positive, real integer values of
%   $s$, in which case we use Series~3,

% \end{itemize}

% The function may be multi-valued for certain $s$. Hence additional
% care must be take along the branch on the real axis for
% $z \in [1,\infty)$, in order to obtain the conventional values (see
% \autoref{sec:branch}). 

% algorithmic2e version
% \begin{algorithm}[t]
%   \KwData{$s, z \in \C$}
%   \KwResult{$Li_s(z)$}
%   $T_1 \leftarrow 0.512$\;
%   $T_2 \leftarrow 10^{-3}$\;
%   $\mu \leftarrow \ln(z)$\; 
%   $t \leftarrow |\mu| /2 \pi$\;
%   \uIf{$2 \pi |z| \leq |\mu|$}{
%     Use Series 1 [defined in \autoref{eq:def}];
%   }
%   \uElseIf{$t \leq T_1$  AND $\Re(s) \leq 0$}{
%     Use Series 2 [defined in \autoref{eq:power_series2a}]\;
%   }
%   \uElseIf{$t \leq T_1$  AND $s$ is further than  $T_2$ from a real integer}{
%     Use Series 2 [defined in \autoref{eq:power_series2a}]\;
%   }
%   \uElseIf{$t \leq T_1$}{
%     Use Series 3 [defined in \autoref{eq:series3} with $j=4$]\;
%   }
%   \Else{
%     Recurse on the Duplication Identity \autoref{eq:duplication}\;
%   }
%   \caption{The polylogarithm algorithm.}
%   \label{alg:main}
% \end{algorithm}

% algorithmicx version
\begin{algorithm}[t]
  \begin{algorithmic}
  % \KwData{$s, z \in \C$}
  % \KwResult{$Li_s(z)$}
  \State $T_1 \leftarrow 0.512$\;
  \State $T_2 \leftarrow 10^{-3}$\;
  \State $\mu \leftarrow \ln(z)$\; 
  \State $t \leftarrow |\mu| /2 \pi$\;
  \If{$2 \pi |z| \leq |\mu|$}
    \State Use Series 1 [defined in \autoref{eq:def}];
  % \ElsIf{$t \leq T_1$  AND $\Re(s) \leq 0$}
  %   \State Use Series 2 [defined in \autoref{eq:power_series2a}]\;
  \ElsIf{$t \leq T_1$  AND $s$ is further than  $T_2$ from a positive, real integer}
    \State Use Series 2 [defined in \autoref{eq:power_series2a}]\;
  \ElsIf{$t \leq T_1$}
    \State Use Series 3 [defined in \autoref{eq:series3} with $j=4$]\;
  \Else
    \State Recurse on the Duplication Identity \autoref{eq:duplication}\;
  \EndIf
  \end{algorithmic}
 \caption{The polylogarithm algorithm $Li_s(z)$ with inputs $s, z \in \C$.}
 \label{alg:main}
\end{algorithm}

\subsection{Stopping Criteria}

Determining at which point to stop each sequence to attain a given
accuracy with minimal cost is not completely trivial. For instance,
Series~1 has terms like $z^k/k^s$. For large $k$ we can approximate
these as $\sim z^k$ and hence the remainder term is that of the
geometric series, and hence for real, positive $z$ it should be almost
trivial to determine a cut-off at which we terminate the series, which
seems to work well for $|z|$ near 0.5. However, for small $z$ the
series drops so fast we don't reach this asymptotic domain. Luckily, a
practical compromise is to terminate the sequence when the relative
value of a summation term drops below 1/2 of the desired precision
bound, \ie at the first $m$ such that
\[
  \frac{z^m/k^m}{\sum_{k=1}^{m} z^k/k^s} \leq 0.5 a, 
\]
where for the majority of the work here we use $a = 10^{-12}$. 
We must be a little careful for very small $z$ and $\Re(s)<0$
to terminate the sequence only after it begins to decrease at
$k = \lceil \Re(s)/\ln |z| \rceil$.

The bound will may be somewhat conservative for imaginary or negative
$z$ where the sequence oscillates and thus may contain cancelling
terms.

We will show the relative merits of this cut-off in the following
sections.

Series~2 appears much more complicated, but has essentially the same
characterisics with the exception that its tail is dominated by powers
of $\ln z$. However, it is somewhat more oscillatory and we have found
that testing the relative size of the last two terms in the tail
are both $\leq 0.5 a$, is more reliable, but once again, it is
somewhat conservative due to the oscillatory nature of the sequence.

The tail of the summation in Series~3 is almost identical to that in
Series~2 and so we use the same termination criteria.

%% file: tests.tex
\section{Tests}

The goal we set here is to attain a relative absolute error
$\leq 10^{-12}$.  This value was chosen to be challenging but
realistic. More importantly, if we set the goal to be machine
precision, we would not be able to see by how much the approach
exceeds the goal and hence is being potentially wasteful of
computation.

% The code allows this to be entered as a parameter, but
% there are a few trade-offs that make automating all aspects of the
% algorithm to achieve a given accuracy with minimal computation
% non-trivial, so we aim for a fixed and useful precision.

The majority of tests performed here were conducted by creating
benchmark data using Mathematica's arbitrary precision {\tt Polylog}
function, to high precision. Note that, in reading this data into
double-precision floating point variables in Julia there is an
inevitable loss of accuracy to machine precision, and so there is a
lower bound on the degree to which we can test the accuracy of our
code, \ie errors around $10^{-15}$ should not be over-interpreted.

In addition to the tests reported below, a large number of additional
tests (nearly 1000) were created using special values, common
identities between the polylogarithm function and other standard
functions, and relationships between polylogarithms, \eg
\autoref{eq:reciprocal_neg_int}.

\subsection{Testing of $|z|$ domains}
\label{sec:test:trade12}

\begin{figure*}[p]
  \centering
  \subfloat[The boundary between Series~1 and Series~2.]{\includegraphics[width=0.495\columnwidth]{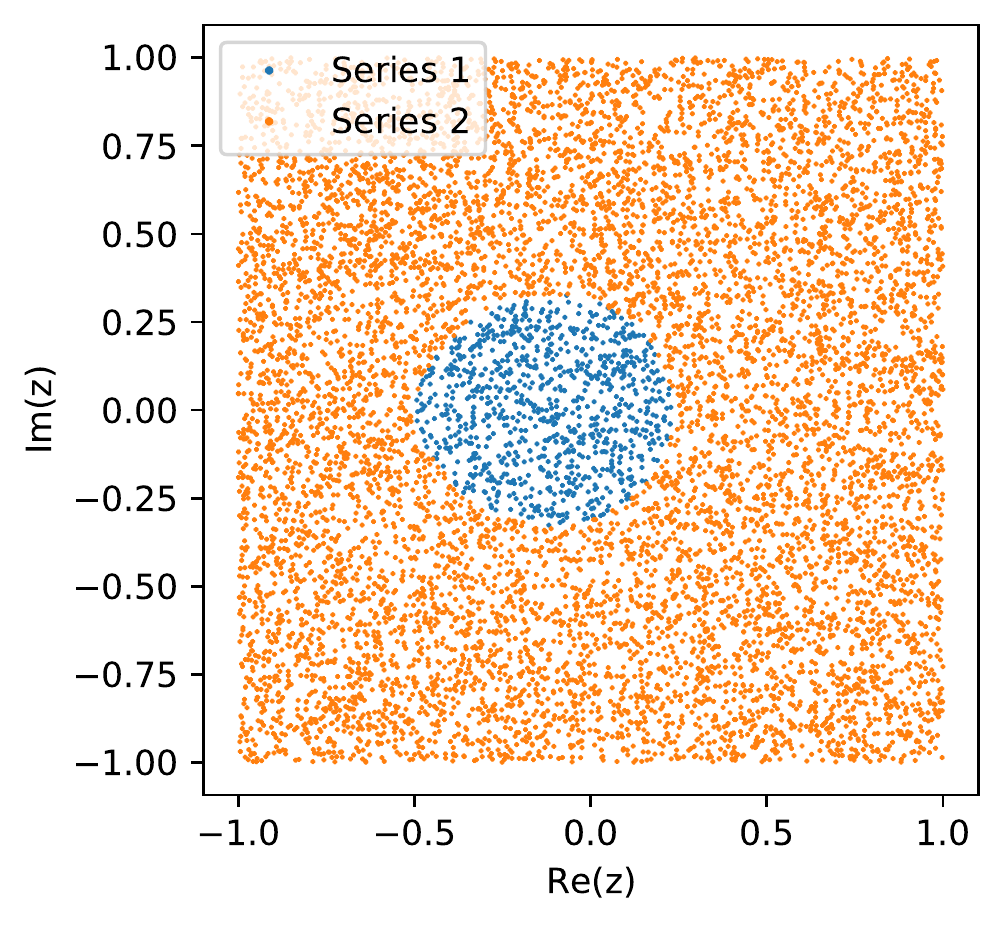}}
  \subfloat[The duplication identity at work.]{\includegraphics[width=0.495\columnwidth]{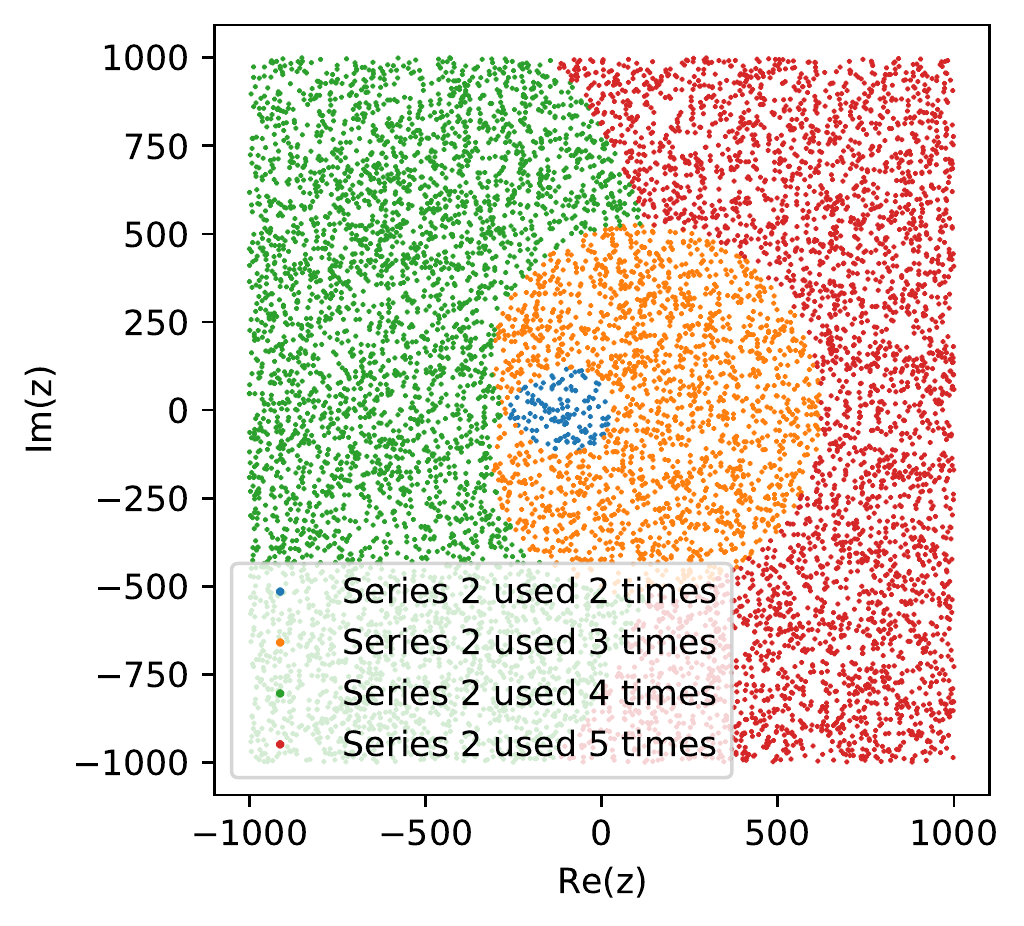}}
  \caption{Illustration of the alternative approaches used as a
    function of $z$.}
  \label{fig:domains}
\end{figure*}

The breakup of the $z$ domain is illustrated in
\autoref{fig:domains}. The left figure shows the boundary (explained
in \autoref{sec:series2}) between use of Series~1 and 2 on a large set
of random points. It is natural to question how well this breakup
works.

A set of test points $z$ were created on the unit disc, spread in even
radii circles about the origin, and we performed this test for around
three different values of $s$ and errors in Series~1 and 2 compared.

\begin{figure*}[p]
  \centering
  % \hspace*{-20mm}%
  \includegraphics[width=0.95\columnwidth]{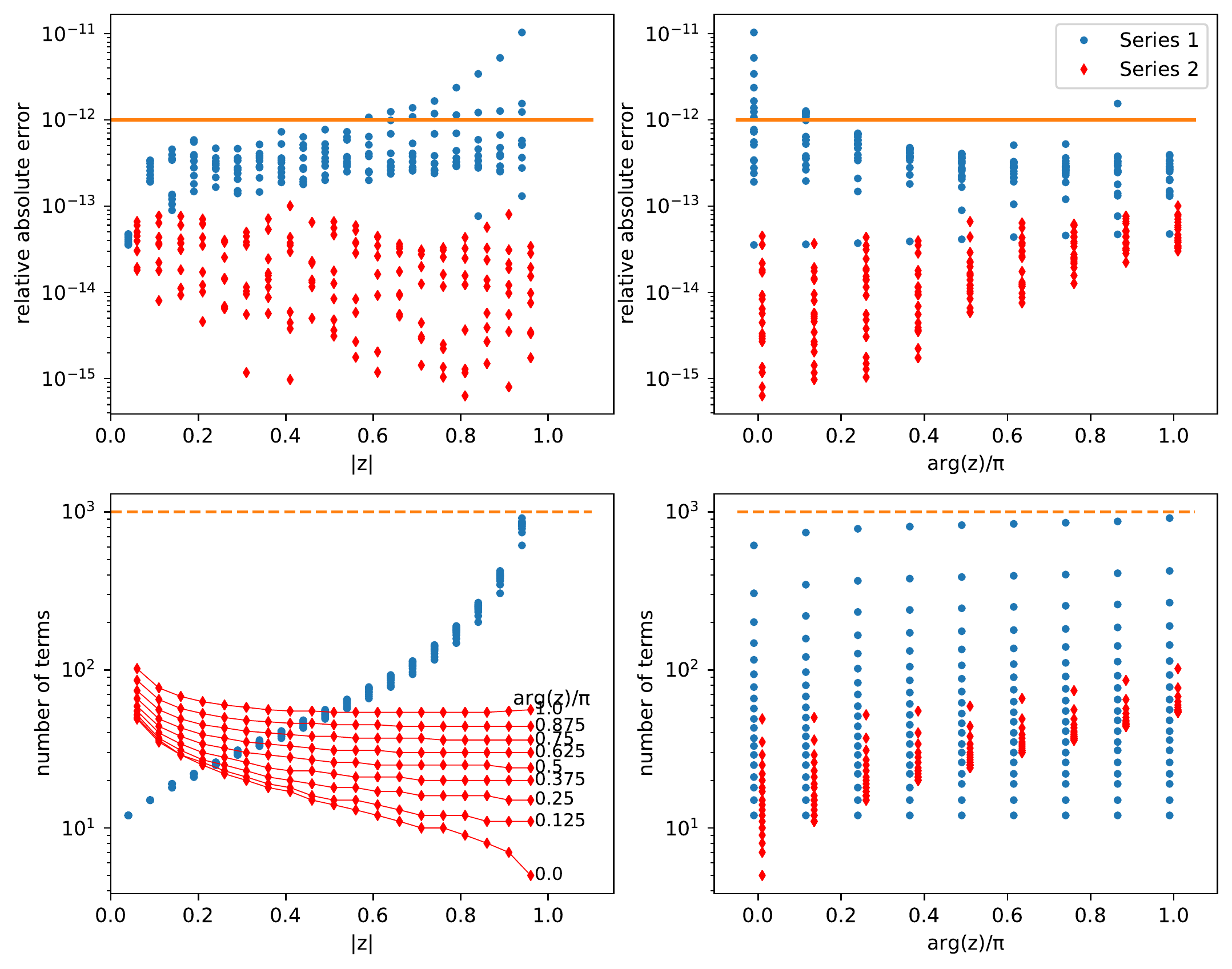}
  \caption{A comparison of Series~1 and 2 for $z$ on the unit circle
    (with $s=-2$). The top two plots show relative accuracy showing
    that the two approaches meet the accuracy requirements within the
    domains in which we calculate them. The stopping rule for Series~2
    is perhaps a little conservative. The two lower plots show number
    of terms calculated in the series showing exactly the trade-off
    described in \autoref{sec:series2}, namely, that the cut-off
    between them lies between approximately 1/4 and 1/2, and that
    while the number of terms in Series~1 is somewhat independent of
    $arg(z)$, Series~2 performs better for positive real values of
    $z$, and worst for negative real values.}
  \label{fig:convergence_a}
\end{figure*}

 \autoref{fig:convergence_a} shows the results: both the relative
absolute errors and number of iterations required by the two series as
functions of $|z|$ and $\arg z$. Noteworthy features include the
trade-off in number of iterations in the two series (lower-left) which
crosses-over between 0.25 and 0.5, with the larger values being along
the negative real axis (bottom right figure).

\begin{figure*}[th]
  \centering
  % \hspace*{-20mm}%
  \includegraphics[width=0.5\columnwidth]{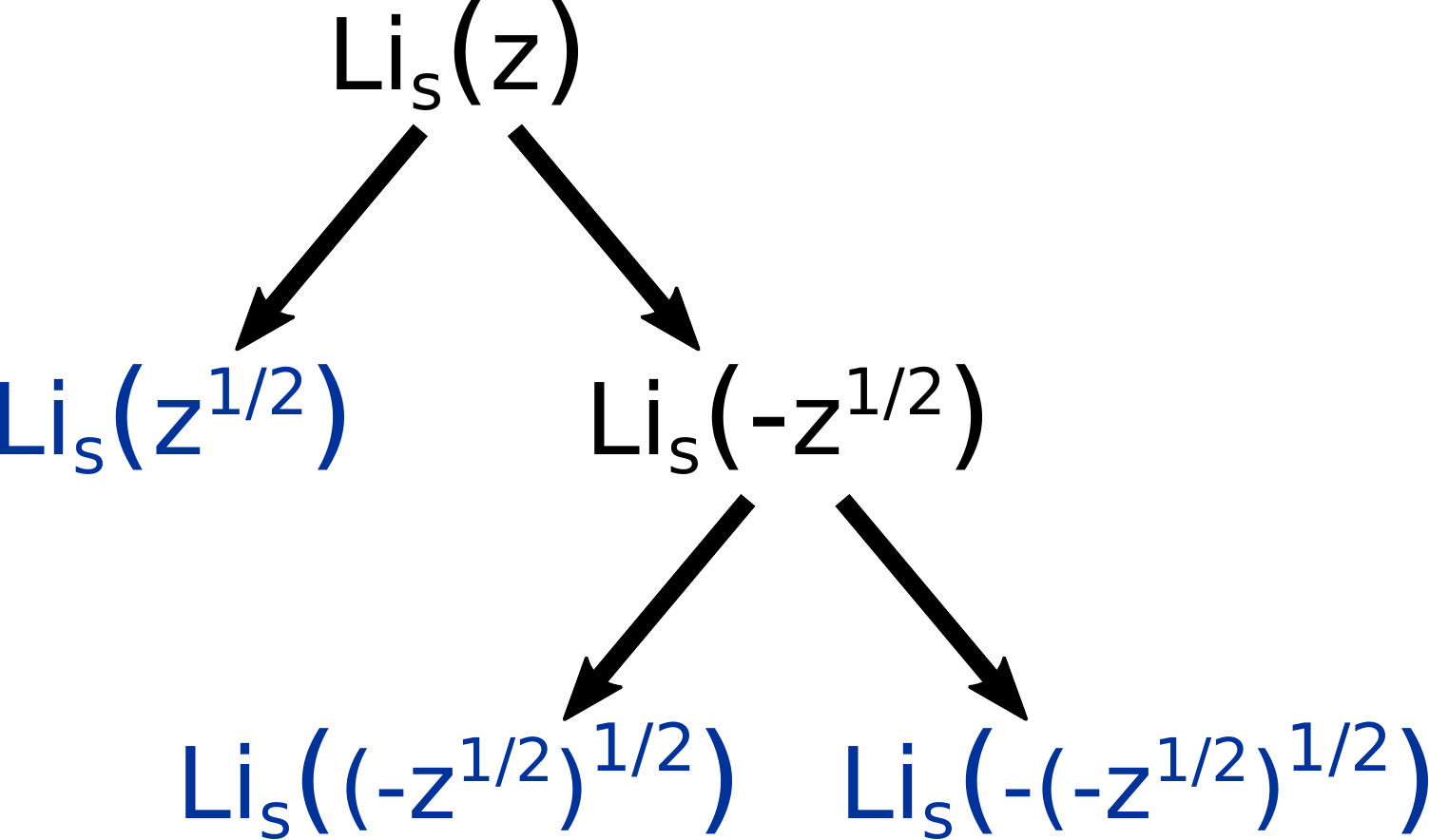}
  \caption{A typical recursion using the duplication identity twice,
    but evaluating Series~2 only 3 times (blue terms), not 4.}
  \label{fig:recursion}
\end{figure*}

\autoref{fig:domains}~(b) illustrates the regions in which the
recursive duplication identity is used, and to what depth. Of note,
the recursion is not a straight-forward matter of breaking the
calculation up into a balanced binary tree to a given depth. The
positive square-root term in the recursion quickly drops into the
domain where we use Series~2. It is the negative square-root term that
bounces around before hitting this region. The recursion tree is
illustrated \autoref{fig:recursion} in a typical case of 2 levels of recursion,
where Series~2 would be evaluated 3 times rather than 4 as one might guess.

\clearpage

\subsection{Accuracy}

We assessed overall accuracy using three sets of 10,000 uniformly
at-random values of $s$ and $z$ with $z$ chosen over 2x2, 16x16 and
2000x2000 rectangles centred on the origin of the complex plane, and
with $s$ chosen in the 16x16 rectangle.  The data used to perform
tests is included with our open source code. The package also contains
code to regenerate all test figures included here.

% $\Re(s), \Im(s) \in [-8,8]$

\begin{figure*}[p]
  \centering
  \subfloat[$z$ in a 2x2 rectangle; 2 points fall outside the goal $10^{-12}$.]{\includegraphics[width=0.65\columnwidth]{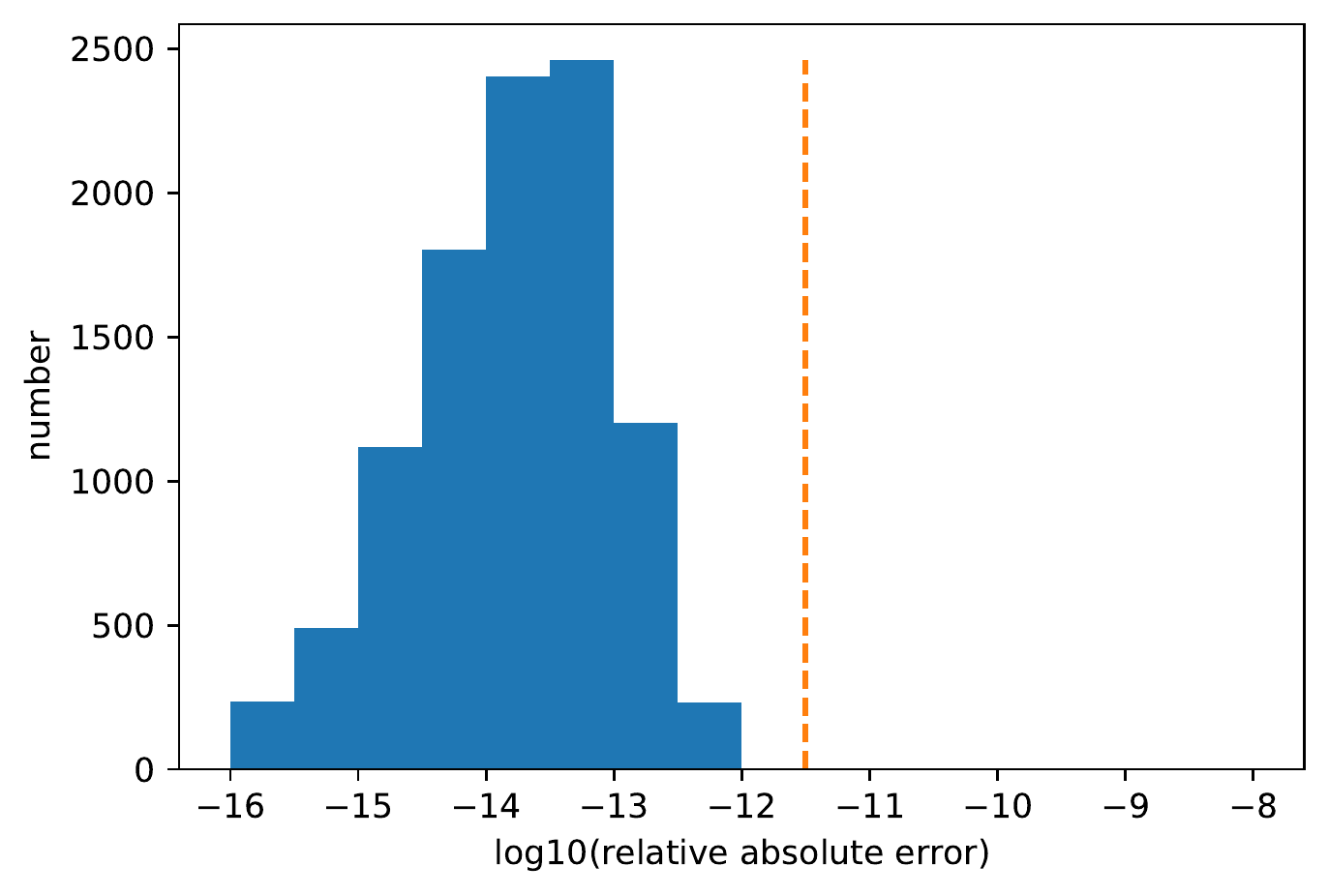}}\\
  \subfloat[$z$ in a 16x16 rectangle; 0 points fall outside the goal $10^{-12}$.]{\includegraphics[width=0.65\columnwidth]{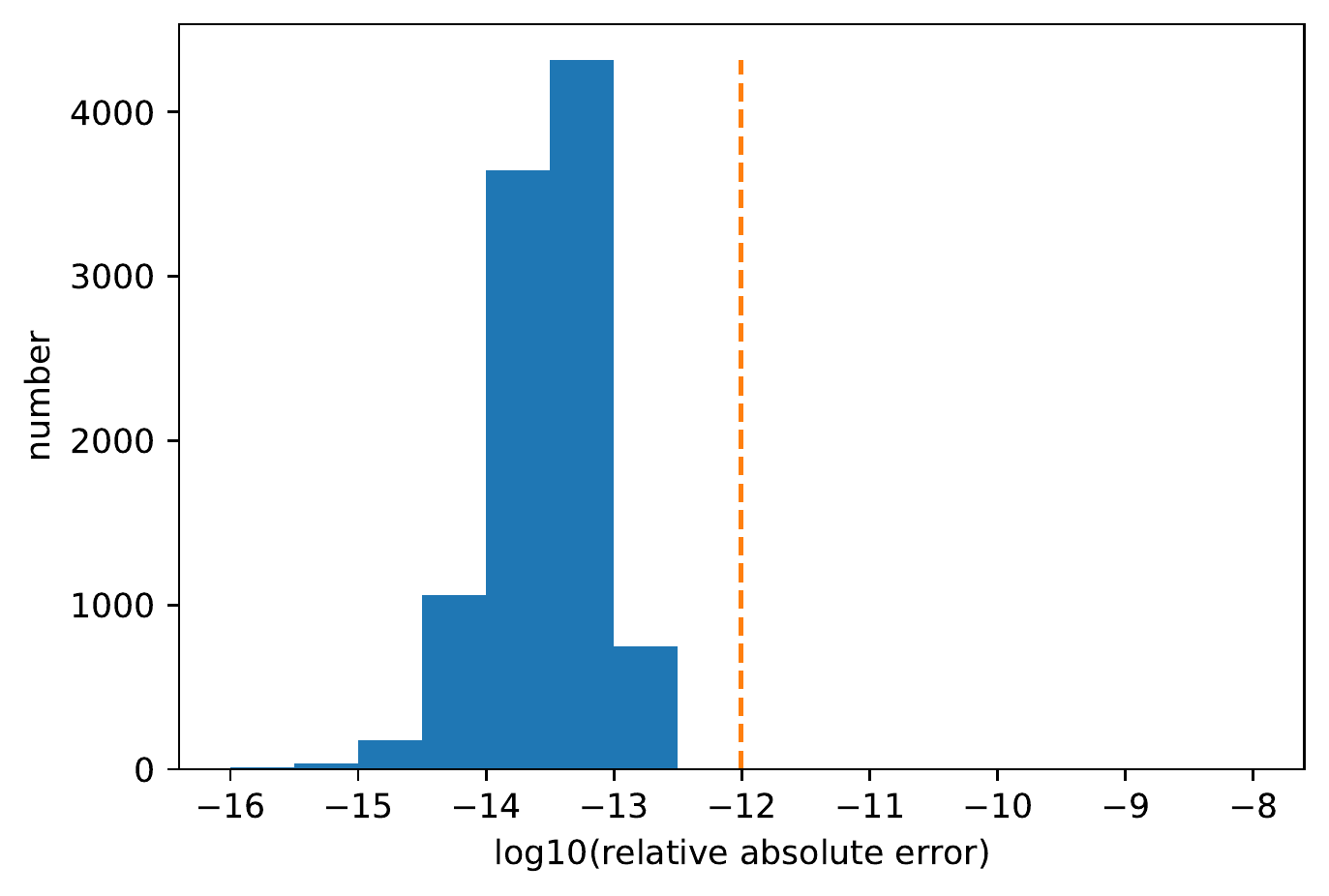}}\\
  \subfloat[$z$ in a 2000x2000 rectangle; 86 points fall outside the goal $10^{-12}$.]{\includegraphics[width=0.65\columnwidth]{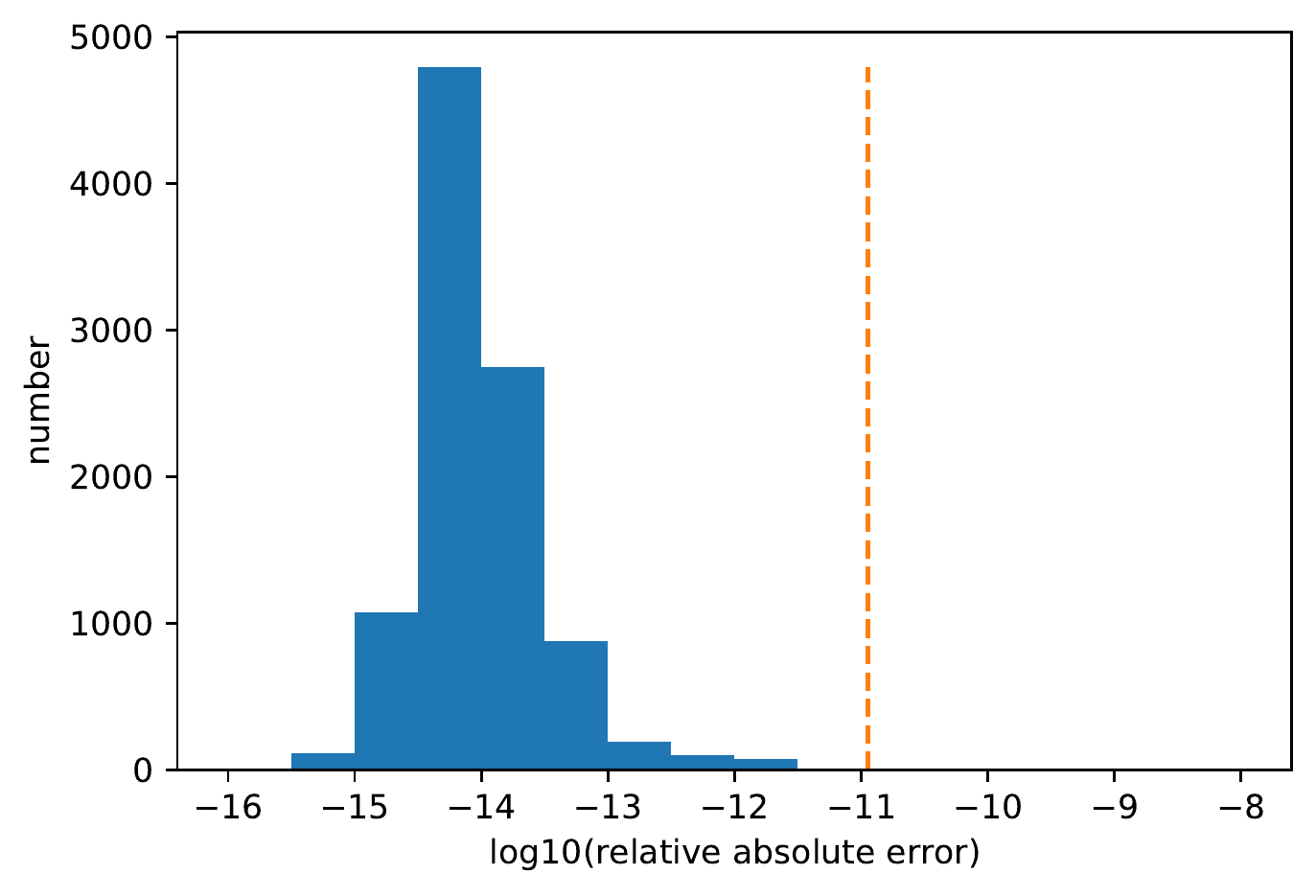}}
  \caption{Histograms of $\log_{10}$ of the relative, absolute
    errors. The dashed lines in each plot indicate the maximum
    error. In total, 88 points out of 30,000 fell outside of the
    desired accuracy goal.}
  \label{fig:accuracy}
\end{figure*}

We present three sets of tests over different sized ranges of $z$ in
order to test the polylogarithm function (i) in the typical region
around the unit circle, (ii) for a larger domain including some
recursion, and (iii) for a very large range of $z$ to stress the
algorithms, particularly the recursion. Histograms of the performance
of the algorithm are shown in \autoref{fig:accuracy}. Only 88 cases
out of 30,000 tests fell outside the desired accuracy of a relative
error no more than $10^{-12}$, and the worst case was
$1.1 \times 10^{-11}$. Achieving a higher degree of accuracy is
possible, but then the majority of points are calculated too
accurately (and hence are wasteful of computations). 

\begin{figure*}[th]
  \centering
  % \hspace*{-20mm}%
  \includegraphics[width=0.5\columnwidth]{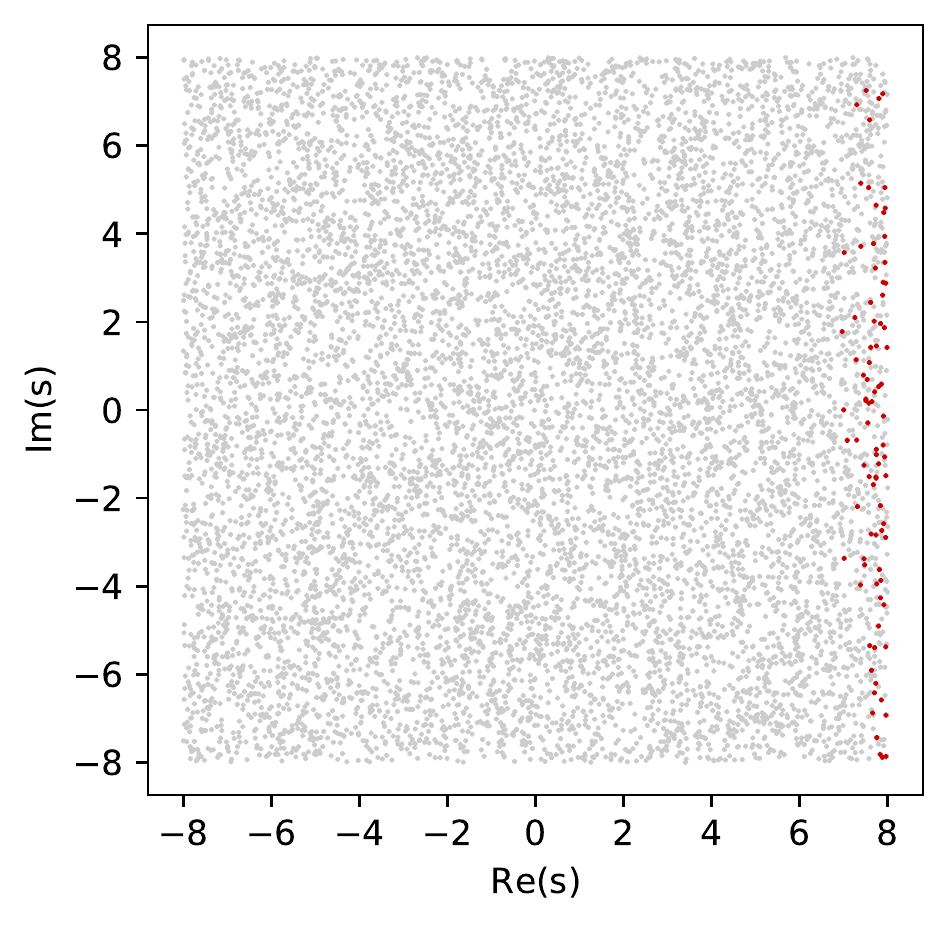}
  \caption{A scatter plot of the random values of $s$ used for
    testing. Red points indicate those with a relative error larger
    than $10^{-12}$.}
  \label{fig:errors_s}
\end{figure*}

It is interesting to understand where the errors occur. The vast
majority occur for large $z$, in the range where 3 or 4 levels of
recursion were required. Moreover, they occur for $s$ with a large,
real component as shown in \autoref{fig:errors_s}.

\subsection{Computational Speed}

\autoref{tab:time} reports the average computation times on the three
random datasets used above on an Intel i9-10900K CPU running in Julia
v1.4.2 with v0.10.3 of the {\tt SpecialFunctions} package, using a
single core, running under Linux Mint 19.3.

The main impact on performance occurs when recursion is needed for
large $z$ values, which requires multiple evaluations of earlier
series.

The table also shows Mathematica\footnote{Calculations for this component
  are made to the same precision goal as our code.} computation
times. The Julia code is at least ten times as fast, which reflects
Julia's desire to provide very fast computations. However, we do not
argue that this is an entirely valid comparison. It is problematic
comparing computation times between programming languages as their
internal measurement facilities may use different concepts and
underlying libraries, and therefore comparisons are not
``apples-to-apples'' comparisons. Nevertheless, the speeds we are
attaining seem creditable for a function of this type.

A second comparison arises in the results presented in
\cite{bailey15:_crandall}, where times for computations of the order
of 20,000-40,000 $\mu$s are reported, though once again note that we
cannot simplistically compare these results as they were computed to
1000 digits of precision. However, presuming 1 bit of accuracy per
summation term, \ie that the computation time is somewhat proportional
to the required accuracy, these times would scale to around 240-480
$\mu$s to obtain 12 digits. The comparison is still flawed because
Bailey and Borwein~\cite{bailey15:_crandall} don't report details of
the system used to perform their computations, but the Julia
implementation seems to compare favourably.

As another comparison point we also report the time to compute Julia's
$\zeta$ function over the same range of inputs. The polylogarithm
calculation is an order of magnitude slower, which is quite reasonable
given that the more common computation strategy (Series~2) requires
evaluation of a series of $\zeta$ functions. However, the result
should allow readers to benchmark the expected performance they would
receive against performance of zeta function calculations in existing
libraries.

\renewcommand{\arraystretch}{1.25}
\begin{table}[t]
  \centering
  \begin{tabular}{r|rrr}
    Dataset & Julia & Mathematica & $\zeta(z)$ \\
    \hline
    $\Re(z), \Im(z) \in [-1,1]$       & 30.3  & 1606.0 & 1.0 \\
    $\Re(z), \Im(z) \in [-8,8]$       & 41.3  & 1790.0 & 1.0 \\
    $\Re(z), \Im(z) \in [-1000,1000]$ & 143.2 & 1890.0 & 0.8 \\
  \end{tabular}
  \caption{Computation times for the polylogaarithm and zeta
    functions. All times are given in $\mu$s and are averages over
    10,000 computations.}
  \label{tab:time}
\end{table}

The current code has room for improvement. The code was written to be
simple and clear and does not exploit Julia's language-dependent
tricks to improve performance. Nor does it use numerical techniques
such as Kahan summation.

The slowest computations involve the recursion. Developing a more
stable reciprocal method would obviate this, and bring these times
back towards those on other domains. However, an interesting
alternative would be to combine the terms in the two components:
Series~2 could then reuse some of the component calculations, notably
the zeta function calculations.

We leave these improvements as future work.

%% file: conclusion.tex
\section{Conclusion}

This paper describes a complete algorithm for computing numerical
values of the polylogarithm for complex arguments using standard
double-precision floating-point calculations.

There are many possibilities for improvements. Most obviously, the
current code has been written for simplicity and clarity and could be
optimised in many ways. Simple improvements include more extensive use
of look-up tables and series termination thresholds that cope better
with oscillation. Moreover series are added in the most direct way as
they are calculated (generally from largest to smallest) which isn't
necessarily optimal, not to mention that Kahan summation or Shank's
transformation might be used to speed up convergence of summations.

There are also additional sequences that could be applied. In
particular, asymptotic forms valid for $|z| \gg 1$, though our
experiments with these showed little promise.

There are many extensions to be made in the future. 
\begin{itemize}
 
\item Extending the toolkit to produce some of the standard related
  functions such as Fermi-Dirac integrals.
  
\item Mass production of calculations: Julia's default approach to
  calculating multiple values of a function is to broadcast the
  inputs, which is ideal for parallelisation, but in calculating
  polylogarithms for multiple values of $z$ but the same $s$, we can
  potentially reuse many of the terms in the calculation, thus saving
  a large amount of computations.
  
\item Calculation of incomplete polylogarithm functions.

\item Calculation of Neilsen generalised polylogarithm
  functions~\cite{jacobs72:_numer_calcul_polyl}. 

\item Calculation of multiple
  polylogarithms~\cite{waldschmidt02:_multip_polyl}. 
 
\end{itemize}